\documentclass[12pt,a4paper]{article}

\usepackage[latin2]{inputenc}
\usepackage{amssymb}
\usepackage{paralist}
\usepackage{amsmath, calc}
\usepackage{amsfonts}
\usepackage{amsthm}
\usepackage{graphicx}
\usepackage{fullpage}
\usepackage{enumerate}
\newtheorem{thm}{Theorem}

\newtheorem{prob}{Problem}

\newtheorem{lem}{Lemma}

\begin{document}

\title{On $B_{h}[1]$-sets which are asymptotic bases of order $2h$}
\author{S\'andor Z. Kiss \thanks{Department of Algebra, Institute of Mathematics, Budapest
University of Technology and Economics, M\H{u}egyetem rkp. 3., H-1111 Budapest, Hungary; ksandor@math.bme.hu.
This author was supported by the National Research, Development and Innovation Office NKFIH Grant No. K129335.
This paper was supported by the J\'anos Bolyai Research Scholarship of the Hungarian Academy of Sciences. Supported by the \'UNKP-20-5 New National Excellence Program
of the Ministry for Innovation and Technology from the source of the National Research Development and Innovation Fund.} , Csaba
S\'andor \thanks{Department of Stohastics, Institute of Mathematics, Budapest
University of Technology and Economics, M\H{u}egyetem rkp. 3., H-1111 Budapest, Hungary; csandor@math.bme.hu.
This author was supported by the Department of Computer Science and Information Theory, Budapest University
of Technology and Economics, M\H{u}egyetem rkp. 3., H-1111 Budapest, Hungary
MTA-BME Lend\"ulet Arithmetic Combinatorics Research Group,
  ELKH, M\H{u}egyetem rkp. 3., H-1111 Budapest, Hungary; This author was supported by the National Research, Development and Innovation Office NKFIH Grant No. K129335.}
}
\date{}
\maketitle

\begin{abstract}
\noindent Let $h,k \ge 2$ be integers. A set $A$ of positive integers is
called asymptotic basis of order
$k$ if every large enough positive integer can be written as the sum of
$k$ terms from $A$. A set of positive integers $A$ is said to be a
$B_{h}[g]$-set if every positive integer can be written as the sum
of $h$ terms from $A$ at most $g$ different ways. In this paper we prove the
existence of $B_{h}[1]$ sets which are asymptotic bases of order $2h$
by using probabilistic methods.


\end{abstract}
\section{Introduction}

Let $h, k \ge 2$ be integers. We denote the set of nonnegative integers by
$\mathbb{N}$ and the set of positive integers by
$\mathbb{Z}^{+}$. Let $A \subset \mathbb{N}$ be an infinite set and let
$R_{h,A}(n)$ denote the number of solutions of the equation
\begin{equation}
a_{1} + a_{2} + \dots + a_{h} = n, \hspace*{3mm} a_{1} \in
A, \dots, a_{h} \in A, \hspace*{3mm} a_{1} \le
a_{2} \le \dots{} \le a_{h},
\end{equation}
\noindent where $n \in \mathbb{N}$. We say a set of positive integers $A$
forms a $B_h[g]$-set if for every $n \in \mathbb{N}$, the number of
representations of $n$ as the sum of $h$ terms in the form (1) is at most $g$,
that is $R_{h,A}(n) \le g$.
A set $A \subset \mathbb{N}$ is said to be an asymptotic basis
of order $k$ if there exists a positive integer $n_{0}$ such that
$R_{k,A}(n) > 0$ for $n > n_{0}$. In \cite{ET} and \cite{ES}, P. Erd\H{o}s,
A. S\'ark\"ozy and V. T. S\'os asked if there exists a Sidon set (i.e., a
$B_2[1]$-set) which is an asymptotic basis of order 3.
It is easy to see that a Sidon set cannot be an asymptotic basis of
order 2 because it does not have enough elements. J. M. Deshouillers and A. Plagne in
\cite{DE} constructed a Sidon set which is an asymptotic basis of order at most
7. In \cite{SB} it was proved the existence of asymptotic
bases of order 5 which are Sidon sets by using probabilistic methods.
In \cite{CE} and \cite{KR} this result was improved on by proving the
existence of an asymptotic basis of order 4 which is a Sidon set. It was also proved \cite{CE} that there exists an asymptotic basis of order 3
which is a $B_2[2]$-set. The above problem of Erd\H{o}s et. al. can be formulated in a more general form.
In a recent paper \cite{KS}, we proved the existence of a $B_{h}[1]$-set which is at the same time forms an asymptotic basis of order $2h+1$.

In this paper we continue the work in this direction.
Particularly, we improve on the above result by proving the
existence of a $B_{h}[1]$-set which is an asymptotic basis of order $2h$.

\begin{thm}
For every $h \ge 2$ integer there exists an asymptotic basis of order $2h$
which is a $B_{h}[1]$-set.
\end{thm}

Before we prove the above theorem, we propose some open problems for further
research. These problems are also appear in \cite{KS}. In general, for $k, h \ge 2$ integers, one can be interested in the existence of an asymptotic basis of
order $k$ which is a $B_{h}[1]$-set. It is easy to see that there does not exist an asymptotic basis of order $k < h$ which is a $B_{h}[1]$-set because it does not have enough elements.
In recent years, it was proved \cite{CRS} that there does not exist a $B_{h}[1]$-set which is an asymptotic basis of order $h$.

\begin{prob}
Determine the smallest value of $k = k(h) > h$ for which there exists a $B_{h}[1]$-set which is an asymptotic basis of order $k$.
\end{prob}
In this paper, we are dealing with
the case $k = 2h$ and we prove the existence of an asymptotic basis of order $2h$  which is a $B_{h}[1]$-set at the same time by using deeper
probabilistic arguments. To prove the existence of an asymptotic basis of order $2h - 1$ which is simultaneously a $B_{h}[1]$-set seems to be very difficult.
In the case when $k = h$, the generalization of the famous conjecture of Erd\H{o}s and Tur\'an states that there does not exist a $B_{h}[g]$-set which is an asymptotic basis of order $h$.
This conjecture is still open and it seems to be hopeless even for $h = 2$.

It is natural question whether there exist an
asymptotic basis of order $h + 1$ which is a $B_{h}[g]$-set for some $g = g(h)$.
For $h \ge 2$, it was proved \cite{SK} the existence of an asymptotic basis of order $h + 1$ which is a $B_{h}[g]$-set.
In \cite{SK}, the order of magnitude of $g = g(h)$ was not controlled. This suggest us to study the following problem.

\begin{prob}
Determine the smallest value of $g = g(h)$ for which there exists an
asymptotic basis of order $h + 1$ which is a $B_{h}[g]$-set.
\end{prob}

In the following section we give a short summary of the probabilistic tools which plays the crucial role in our proof.

\section{Probabilistic tools}

In the proof of Theorem 1, we apply the probabilistic method due to Erd\H{o}s and R\'enyi. Interested reader can find a nice summary about this method in
the book of Halberstam and Roth \cite{HL}.
Let $\Omega$ denote the set of strictly increasing sequences of
positive integers. In this paper we denote the probability of an event
$\mathcal{E}$ by $\mathbb{P}(\mathcal{E})$, and the expectation of a random
variable $\zeta$ by $\mathbb{E}(\zeta)$.

\begin{lem}
Let
\[
\alpha_{1}, \alpha_{2}, \alpha_{3} \dots{}
\]
be real numbers satisfying
\[
0 \le \alpha_{n} \le 1 \hspace*{4mm} (n = 1, 2, \dots{}).
\]
Then there exists a probability space ($\Omega$, $\mathcal{X}$, $\mathbb{P}$) with the following two properties:
\begin{itemize}
\item[(i)] For every natural number $n$, the event $\mathcal{E}^{(n)} = \{\mathcal{A}$:
  $\mathcal{A} \in \Omega$, $n \in \mathcal{A}\}$ is measurable, and
  $\mathbb{P}(\mathcal{E}^{(n)}) = \alpha_{n}$.
\item[(ii)] The events $\mathcal{E}^{(1)}$, $\mathcal{E}^{(2)}$, ... are independent.
\end{itemize}
\end{lem}
See Theorem 13. in \cite{HL}, p. 142.

We denote the characteristic function of the event $\mathcal{E}^{(n)}$ by $t_{n}$ or we can say the boolean random variable means that:
\[
t_{n} =
\left\{
\begin{aligned}
1 \textnormal{, if } n \in \mathcal{A} \\
0 \textnormal{, if } n \notin \mathcal{A}
\end{aligned} \hspace*{3mm}.
\right.
\]

\noindent Moreover, we denote the number of solutions of the equation
$a_{i_{1}} + a_{i_{2}} + \dots{} + a_{i_{k}} = n$ by $r_{k,\mathcal{A}}(n)$,
where
$a_{i_{1}} \in \mathcal{A}$, $a_{i_{2}} \in \mathcal{A}$, ...,$a_{i_{k}} \in \mathcal{A}$, $1 \le a_{i_{1}} < a_{i_{2}} \dots{}
< a_{i_{k}} < n$.
Thus we have
\[
r_{k,\mathcal{A}}(n) = \sum_{\overset{(a_{1}, a_{2}, \dots{}, a_{k}) \in
\mathbb{N}^{k}}{1
    \le a_{1} < \dots{} < a_{k} < n}\atop {a_{1} + a_{2} + \dots{} + a_{k} =
    n}}t_{a_{1}}t_{a_{2}} \dots{} t_{a_{k}}.
\]
It is easy to see that $r_{k,\mathcal{A}}(n)$ is the sum of random
variables. It is clear that for $k > 2$ these variables are not necessarily independent because the
same $t_{a_{i}}$ may appear in many terms. To handle this
problem we need more advanced probabilistic tools.

Our proof is based on a method of J. H. Kim and V. H. Vu \cite{JK}, \cite{TA}, \cite{VA}, \cite{VU}.
Assume that $t_{1}, t_{2}, \dots{} , t_{n}$ are
independent binary (i.e., all $t_{i}$'s are in $\{0,1\}$) random variables.
Consider a polynomial $Y = Y(t_{1}, t_{2}, \dots{} ,
t_{n})$ in $t_{1}, t_{2}, \dots{} ,
t_{n}$ with degree $k$. A polynomial $Y$ is said to be
totally positive if it can be
written in the form $Y = \sum_{i}e_{i}\Gamma_{i}$, where the $e_{i}$'s are
positive and $\Gamma_{i}$ is a product of some $t_{j}$'s.
Given any multi-index
$\overline{\eta} = (\eta_{1}, \dots{} ,\eta_{n}) \in \mathbb{N}^{n}$,
we define the partial derivative $\partial^{(\overline{\eta})}(Y)$ of
$Y$ by
\[
\partial^{(\overline{\eta})}(Y) = \left(\frac{\partial}{\partial t_{1}}\right)^{\eta_{1}} \cdots{} \left(\frac{\partial}{\partial t_{n}}\right)^{\eta_{n}}Y(t_{1}, t_{2}, \dots{} ,t_{n}),
\]
and denote the order of $\overline{\eta}$ as $|\overline{\eta}| = \eta_{1} + \dots{} + \eta_{n}$. For any order $d \ge 0$, we denote $\mathbb{E}_{d}(Y) = \max_{\eta: |\eta| = d}\mathbb{E}(\partial^{(\overline{\eta})}(Y))$. Then $\mathbb{E}_{0}(Y) = \mathbb{E}(Y)$ and $\mathbb{E}_{d}(Y) = 0$ if $d$ exceeds the
degree of $Y$. Define $\mathbb{E}_{\ge d}(Y) =  \max_{d^{'} \ge d}\mathbb{E}_{d^{'}}(Y)$. We will apply the following theorem proved by Kim and Vu, which informally states that when the all the partial derivatives of a totally positive polynomial $Y$ of degree $k$ are less on average than $Y$ itself and $k$ is small in some sense, then $Y$ is concentrated around its mean.

\begin{lem}(J. H. Kim - V. H. Vu)
For every positive integer $k$ and $Y = Y(t_{1}, t_{2}, \dots{} ,t_{n})$
totally positive polynomial of degree $k$, where
the $t_{i}$'s are independent binary random variables, and for any
$\lambda > 0$ there exists a constant $d_{k} > 0$ depending only on $k$ such
that
\[
\mathbb{P}\left(|Y-\mathbb{E}(Y)| \ge d_{k}\lambda^{k-\frac{1}{2}}\sqrt{\mathbb{E}_{\ge 0}(Y)\mathbb{E}_{\ge 1}(Y)}\right) = O_{k}\left(e^{-\lambda/4+(k-1)\log n}\right).
\]
\end{lem}

\noindent See \cite{VU} for the proof. Finally, we need the Borel - Cantelli lemma which is Theorem 7. in \cite{HL}, p. 135.

\begin{lem}(Borel - Cantelli)
Let $\{B_{i}\}$ be a sequence of events in a probability space. If
\[
\sum_{j=1}^{+\infty}\mathbb{P}(B_{j}) < \infty,
\]
\noindent then with probability 1, at most a finite number of the events
$B_{j}$ can occur.
\end{lem}
\noindent In the next section we prove a lemma which plays a very important role in the proof of Theorem 1.

\section{An auxiliary tool}
Throughout the remainig part of the
paper we use the notation $f(x) \ll g(x)$ which means $f(x) = O(g(x))$. We also use the notation $f(x) \asymp g(x)$ which means $f(x) = \Theta(g(x))$. Next, we prove the following technical lemma which plays
an important role in the proofs.

\begin{lem}
\begin{itemize}
\item[(i)] Let $M \ge 2$ be a positive integer, and let
$0 < \alpha, \beta < 1$ be arbitrary real numbers. Then
\[
\sum_{n=1}^{M-1}\frac{1}{n^{\alpha}}\cdot \frac{1}{(M - n)^{\beta}} \ll
\frac{1}{M^{\alpha+\beta-1}}.
\]
\item[(ii)] Let $M$ be an arbitrary integer, and let $0 < \alpha, \beta < 1$
with $\alpha + \beta > 1$ be arbitrary real numbers. Then
\[
\sum_{n=1}^{\infty}\frac{1}{(|n + M| + 1)^{\alpha}}\cdot \frac{1}{n^{\beta}}
\ll \frac{1}{(|M|+1)^{\alpha+\beta-1}}.
\]
\item[(iii)] Let $l \le 2h$. Then, for every positive integer $M$,
\[
\sum_{\substack{(z_{1},\dots{} ,z_{l})\in (\mathbb{Z}^{+})^{l} \\ z_{1} + \dots{} + z_{l} = M}}
\frac{1}{(z_{1}\cdots{} z_{l})^{\frac{4h-3}{4h-1}}}\ll \frac{1}{M^{1-\frac{2l}{4h-1}}}.
\]
\item[(iv)] Let $0 \le s \le t \le 2h$. Then, for every integer $M$,
\[
\sum_{\substack{(z_{1},\dots{} ,z_{t})\in (\mathbb{Z}^{+})^{t} \\ z_{1} + \dots{} + z_{s} - (z_{s+1} + \dots{} + z_{t}) = M}}
\frac{1}{(z_{1}\cdots{} z_{t})^{\frac{4h-3}{4h-1}}}\ll \frac{1}{(|M|+1)^{1-\frac{2t}{4h-1}}}.
\]
\end{itemize}
\end{lem}

\subsection{Proof of Lemma 4}
We estimate the sums by integral according to the well known
Euler integral formula.
\[
\sum_{n=1}^{M-1}\frac{1}{n^{\alpha}}\cdot \frac{1}{(M - n)^{\beta}} = \sum_{n=1}^{\lfloor M/2 \rfloor}\frac{1}{n^{\alpha}}\cdot \frac{1}{(M - n)^{\beta}} + \sum_{n=\lfloor M/2 \rfloor + 1}^{M-1}\frac{1}{n^{\alpha}}\cdot \frac{1}{(M - n)^{\beta}}
\]
\[
\ll \frac{1}{M^{\beta}}\sum_{n=1}^{\lfloor M/2 \rfloor}\frac{1}{n^{\alpha}} + \frac{1}{M^{\alpha}}\cdot \sum_{n=\lfloor M/2 \rfloor + 1}^{M-1}\frac{1}{(M - n)^{\beta}}
\ll \frac{1}{M^{\beta}}\int_{0}^{M}\frac{dx}{x^{\alpha}} + \frac{1}{M^{\alpha}}\int_{0}^{M}\frac{dx}{x^{\beta}} \ll \frac{1}{M^{\alpha+\beta-1}},
\]
which proves (i).

Moreover, if $M \ge 0$,
\[
\sum_{n=1}^{\infty}\frac{1}{(|n + M| + 1)^{\alpha}}\cdot \frac{1}{n^{\beta}} = \sum_{n=1}^{M}\frac{1}{(n + M + 1)^{\alpha}}\cdot \frac{1}{n^{\beta}} + \sum_{n=M+1}^{\infty}\frac{1}{(n + M + 1)^{\alpha}}\cdot \frac{1}{n^{\beta}}
\]
\[
\le \frac{1}{(M+1)^{\alpha}}\sum_{n=1}^{M}\frac{1}{n^{\beta}} + \sum_{n=M+1}^{\infty}\frac{1}{n^{\alpha+\beta}}
\]
\[
\ll \frac{1}{(M+1)^{\alpha}}\int_{0}^{M+1}\frac{dx}{x^{\beta}} + \int_{M+1}^{\infty}\frac{dx}{x^{\alpha+\beta}} \ll \frac{1}{(|M|+1)^{\alpha+\beta-1}}.
\]
Furthermore, if $M < 0$,
\[
\sum_{n=1}^{\infty}\frac{1}{(|n + M| + 1)^{\alpha}}\cdot \frac{1}{n^{\beta}} = \sum_{n=1}^{-\lfloor M/2 \rfloor}\frac{1}{(|M + n| + 1)^{\alpha}}\cdot \frac{1}{n^{\beta}} + \sum_{n=-\lfloor M/2 \rfloor + 1}^{-M}\frac{1}{(|M + n| + 1)^{\alpha}}\cdot \frac{1}{n^{\beta}}
\]
\[
+ \sum_{n=-M+1}^{-2M}\frac{1}{(|M + n| + 1)^{\alpha}}\cdot \frac{1}{n^{\beta}}
+ \sum_{n=-2M+1}^{\infty}\frac{1}{(|M + n| + 1)^{\alpha}}\cdot \frac{1}{n^{\beta}}
\]
\[
\ll \frac{1}{|M|^{\alpha}}\sum_{n=1}^{-\lfloor M/2 \rfloor}\frac{1}{n^{\beta}} + \frac{1}{|M|^{\beta}}\sum_{n=-\lfloor M/2 \rfloor + 1}^{-M}\frac{1}{(|M + n| + 1)^{\alpha}} + \frac{1}{|M|^{\beta}}\sum_{n=-M+1}^{-2M}\frac{1}{(|M + n| + 1)^{\alpha}}
\]
\[
+ \sum_{n=-M+1}^{\infty}\frac{1}{n^{\alpha+\beta}} \ll \frac{1}{|M|^{\alpha}}\int_{0}^{-M}\frac{dx}{x^{\beta}} + \frac{1}{|M|^{\beta}}\int_{0}^{-M}\frac{dx}{x^{\alpha}} + \frac{1}{|M|^{\beta}}\int_{0}^{-M}\frac{dx}{x^{\alpha}} + \int_{-M}^{\infty}\frac{dx}{x^{\alpha+\beta}}
\]
\[
\ll \frac{1}{(|M| + 1)^{\alpha+\beta-1}},
\]
which proves (ii).

Next, we prove (iii) by induction on $l$. The statement is clear for $l = 1$. Assume that it is true for $l - 1 \le 2h - 1$. Then by (i),
\[
\sum_{\substack{(z_{1},\dots{} ,z_{l})\in (\mathbb{Z}^{+})^{l} \\ z_{1} + \dots{} + z_{l} = M}}
\frac{1}{(z_{1}\cdots{} z_{l})^{\frac{4h-3}{4h-1}}} = \sum_{z_{l}=1}^{M-1}\frac{1}{z_{l}^{\frac{4h-3}{4h-1}}}\sum_{\substack{(z_{1},\dots{} ,z_{l-1})\in (\mathbb{Z}^{+})^{l-1} \\ z_{1} + \dots{} + z_{l-1} = M-z_{l}}}
\frac{1}{(z_{1}\cdots{} z_{l-1})^{\frac{4h-3}{4h-1}}}
\]
\[
\ll \sum_{z_{l}=1}^{M-1}\frac{1}{z_{l}^{\frac{4h-3}{4h-1}}}\cdot \frac{1}{(M - z_{l})^{1-\frac{2(l-1)}{4h-1}}} \ll \frac{1}{M^{1-\frac{2l}{4h-1}}},
\]
which gives (iii).

Finally, if either $s = 0$, or $t = s$, the statement of (iv) follows immediately from (iii). Otherwise, one can assume that $0 < s < t$. It follows from (ii) and (iii) that
\[
\sum_{\substack{(z_{1},\dots{} ,z_{t})\in (\mathbb{Z}^{+})^{t} \\ z_{1} + \dots{} + z_{s} - (z_{s+1} + \dots{} + z_{t}) = M}}\frac{1}{(z_{1}\cdots{} z_{t})^{\frac{4h-3}{4h-1}}}
\]
\[
= \sum_{n=1}^{\infty}\left(\sum_{\substack{(z_{1},\dots{} ,z_{s})\in (\mathbb{Z}^{+})^{s} \\ z_{1} + \dots{} + z_{s} = n + M}}\left(\frac{1}{(z_{1}\cdots{} z_{s})^{\frac{4h-3}{4h-1}}}\sum_{\substack{(z_{s+1},\dots{} ,z_{t})\in (\mathbb{Z}^{+})^{t-s} \\ z_{s+1} + \dots{} + z_{t} = n}}\frac{1}{(z_{s+1}\cdots{} z_{t})^{\frac{4h-3}{4h-1}}}\right)\right)
\]
\[
\ll \sum_{n=1}^{\infty}\sum_{\substack{(z_{1},\dots{} ,z_{s})\in (\mathbb{Z}^{+})^{s} \\ z_{1} + \dots{} + z_{s} = n + M}}\frac{1}{(z_{1}\cdots{} z_{s})^{\frac{4h-3}{4h-1}}}\cdot \frac{1}{n^{1-\frac{2(t-s)}{4h-1}}}
\ll \sum_{n=1}^{\infty}\frac{1}{(|n+M|+1)^{1-\frac{2s}{4h-1}}}\cdot \frac{1}{n^{1-\frac{2(t-s)}{4h-1}}}
\]
\[
\ll \frac{1}{(|M|+1)^{1-\frac{2t}{4h-1}}},
\]
which proves (iv). The proof of Lemma 4 is completed.

\section{Proof of Theorem 1}

\subsection{Outline of the proof}

Let $h$ be fixed and let $\alpha = \frac{2}{4h-1}$. Define the sequence
$\alpha_{n}$ in Lemma 1 by

\[
\alpha_{n} = \frac{1}{n^{1-\alpha}},
\]

\noindent so that $\mathbb{P}(\{B$:
  $B \in \Omega$, $n \in B\}) =
  \frac{1}{n^{1-\alpha}}$.
We prove that in this probability space, almost always one can remove infinitely many elements from a set $B\in \Omega$ such that the remaining set
is both a $B_{h}[1]$ set and an asymptotic basis of order $2h$. In the first step, we show that almost surely $R_{2h,B}(n)$ is large, if $n$ is large enough.
Next, we delete elements from the set $B$ to get a $B_{h}[1]$-set $A$. Finally, we prove that with probability 1, $R_{2h,A}(n)$ is still large, if $n$ is large enough, which implies that the set $A$ is suitable.

We prove that $B$ is an asymptotic basis of order $2h$ in the following stronger sense. There exists a constant $C_{h} > 0$ such that

\begin{equation}
R_{2h,B}(n) > C_{h}n^{\frac{1}{4h-1}}
\end{equation}

for every $n$ large enough, with probability 1.

To do this, we apply the following lemma proved in \cite{SK}
with $S = A$, $k = 2h$ and $\alpha = \frac{2}{4h-1}$.

\begin{lem}
Let $k \ge 2$ be a fixed integer and let $\mathbb{P}(\{A$:
  $A \in \Omega$, $n \in A\}) =
\frac{1}{n^{1-\alpha}}$ where $\alpha > \frac{1}{k}$. Then with probability 1,
 $r_{k,A}(n) > cn^{k\alpha-1}$ for every sufficiently large $n$,
where $c = c(\alpha,k)$ is a suitable positive constant.
\end{lem}

Obviously, $R_{2h,B}(n) \ge r_{2h,B}(n)$ with probability 1, so the proof of (2) is completed. To get a $B_{h}[1]$ set, one has to remove the elements which hurt the $B_{h}[1]$ property. The $B_{h}[1]$ property can be hurt in two different ways. First, if there exists two $h$-tuples formed by pairwise distinct elements of $B$ such that the sums of the terms in both $h$-tuples are the same i.e.,
there exist pairwise distinct $b_{1},\dots{} ,b_{2h}\in B$ with
\[
b_{1} + \dots{} + b_{h} = b_{h+1} + \dots{} + b_{2h}.
\]
On the other hand, if there exist $b_{1},\dots{} ,b_{2h}\in B$ with $b_{1} + \dots{} + b_{h} = b_{h+1} + \dots{} + b_{2h}$, where $b_{1},\dots{} ,b_{2h}$ are not distinct,
then by subtracting the common terms from both sides and collecting the equal terms, the previous equality become of the form
\[
d_{1}b_{1}^{'} + d_{2}b_{2}^{'} + \dots{} + d_{k}b_{k}^{'} = e_{1}b_{k+1}^{'} + \dots{} + e_{l}b_{k+l}^{'},
\]
with positive integer weights $d_{1},\dots{} ,d_{k}, e_{1},\dots{} ,e_{l}\in \mathbb{Z}^{+}$, $k + l \le 2h - 1$, $d_{1} + \dots{} + d_{k} = e_{1} + \dots{} + e_{l} \le h$
and $b_{1}^{'},\dots{} ,b_{k+l}^{'}\in B$ are already pairwise distinct.

We will prove that by removing the largest element from the above equalities, then with probability 1, the remaining set will be suitable. More
formally, define the set $C$ by
\[
C = \{b: b\in B, \exists b_{2},\dots{} ,b_{2h}\in B \textnormal{ distinct }, b_{i} < b, \hspace*{1.5mm} b + b_{2} + \dots{} + b_{h} = b_{h+1} + \dots{} + b_{2h}\}
\]
\[
\cup \{b: b\in B, \exists d_{1},\dots{} ,d_{k}, e_{1},\dots{} ,e_{l} \in \mathbb{Z}^{+}, k + l \le 2h - 1, d_{1} + \dots{} + d_{k} = e_{1} + \dots{} + e_{l}\le h
\]
\[
\exists b_{2},\dots{} ,b_{k+l}\in B \textnormal{ distinct }, b_{i} < b, \hspace*{1.5mm} d_{1}b + d_{2}b_{2} + \dots{} + d_{k}b_{k} = e_{1}b_{k+1} + \dots{} + e_{l}b_{k+l}\}.
\]
We show that $A = B\setminus C$ is suitable with probability 1. To prove this,
we need to show that there exists an $n_{2}\in \mathbb{Z}^{+}$ such that
\begin{equation}
R_{2h,A}(n) \ge 1
\end{equation}
for every $n \ge n_{2}$, with probability 1. We write
\begin{equation}
R_{2h,A}(n) = R_{2h,B\setminus C}(n) = R_{2h,B}(n) - (R_{2h,B}(n) - R_{2h,B\setminus C}(n)).
\end{equation}
Thus we need an upper estimation to $R_{2h,B}(n) - R_{2h,B\setminus C}(n)$. We have
\[
R_{2h,B}(n) - R_{2h,B\setminus C}(n) = |\{(b_{1},\dots{} ,b_{2h}): b_{1},\dots{} ,b_{2h}\in B, b_{1}\le \dots{} \le b_{2h},
b_{1} + \dots{} + b_{2h} = n,
\]
\[
\exists 1 \le i \le 2h, b_{i}\in C\}|
\]
\[
\le |\{(b_{1},\dots{} ,b_{2h}): b_{1},\dots{} ,b_{2h}\in B, b_{1}\le \dots{} \le b_{2h},
b_{1} + \dots{} + b_{2h} = n, \textnormal{ there are equal terms among }
\]
\[
b_{1},\dots{} ,b_{2h}\}|
\]
\[
+ |\{(b_{1},\dots{} ,b_{2h}): b_{1},\dots{} ,b_{2h}\in B, b_{1} < \dots{} < b_{2h},
b_{1} + \dots{} + b_{2h} = n, \exists 1\le i \le 2h, \exists b_{2}^{'},\dots{} ,b_{2h}^{'}\in B
\]
\[
\textnormal{ distinct }, b_{j}^{'} < b_{i}, b_{i} + b_{2}^{'} + \dots{} + b_{h}^{'} = b_{h+1}^{'} + \dots{} + b_{2h}^{'}\}|
\]
\[
+ |\{(b_{1},\dots{} ,b_{2h}): b_{1},\dots{} ,b_{2h}\in B, b_{1} < \dots{} < b_{2h},
b_{1} + \dots{} + b_{2h} = n, \exists 1\le i \le 2h,
\]
\[
\exists d_{1},\dots{} ,d_{k}, e_{1},\dots{} ,e_{l} \in \mathbb{Z}^{+}, k + l \le 2h - 1, d_{1} + \dots{} + d_{k} = e_{1} + \dots{} + e_{l}\le h,
\]
\[
\exists b_{2}^{'},\dots{} ,b_{k+l}^{'}\in B \textnormal{ pairwise distinct }, b_{j}^{'} < b_{i}, d_{1}b_{i} + d_{2}b_{2}^{'} + \dots{} + d_{k}b_{k}^{'} = e_{1}b_{k+1}^{'} + \dots{} + e_{l}b_{k+l}^{'}\}|
\]
\[
= R_{B}^{(1)}(n) + R_{B}^{(2)}(n) + R_{B}^{(3)}(n).
\]

\noindent \textbf{Case 1.} Upper estimation for $R_{B}^{(1)}(n)$. It is clear that
\[
R_{B}^{(1)}(n) = \sum_{\substack{(f_{1}, \dots{} ,f_{t})\in (\mathbb{Z}^{+})^{t}\\ f_{1} + \dots{} + f_{t} = 2h\\ t \le 2h - 1}}|\{(b_{1},\dots{} ,b_{t}): b_{1},\dots{} ,b_{t}\in B, b_{1} < \dots{} < b_{t}, f_{1}b_{1} + \dots{} + f_{t}b_{t} = n\}|.
\]
The following lemma ensures that with probability 1, there are only bounded number of such representations.

\begin{lem}
For positive integers $t \le 2h - 1$ and $f_{1}, \dots{} ,f_{t}\in \mathbb{Z}^{+}$ with
$f_{1} + \dots{} + f_{t} \le 2h$, the number of solutions of the equation
$f_{1}x_{1} + \dots{} + f_{t}x_{t} = m$, where $x_{1}, \dots{} ,x_{t} \in B$
are distinct, is almost always bounded for every positive integer $m$.
\end{lem}

It follows from Lemma 7 that with probability 1, there exists a positive integer $C_{f_{1}, \dots{} ,f_{t}}(B)$ depending on the set $B$ such that
\[
R_{B}^{(1)}(n) \le \sum_{\substack{(f_{1}, \dots{} ,f_{t})\in (\mathbb{Z}^{+})^{t}\\ f_{1} + \dots{} + f_{t} = 2h\\ t \le 2h - 1}}C_{f_{1}, \dots{} ,f_{t}}(B) = C^{(1)}(B).
\]

\noindent \textbf{Case 2.} Upper estimation for $R_{B}^{(2)}(n)$. Obviously,
\[
R_{B}^{(2)}(n) = \sum_{u=1}^{h}\sum_{v=u}^{u+h}|\{(b_{1},\dots{} ,b_{2h}): b_{1},\dots{} ,b_{2h}\in B, b_{1} < \dots{} < b_{2h},
b_{1} + \dots{} + b_{2h} = n,
\]
\[
\exists 1\le i \le 2h, \exists b_{2}^{'},\dots{} ,b_{2h}^{'}\in B \textnormal{ distinct }, b_{j}^{'} < b_{i}, 2 \le j \le 2h,
\]
\[
b_{i} + b_{2}^{'} + \dots{} + b_{h}^{'} = b_{h+1}^{'} + \dots{} + b_{2h}^{'}, u = |\{b_{i}, b^{'}_{2}, \dots{} ,b^{'}_{h}\}\cap \{b_{1}, \dots{} ,b_{2h}\}|,
\]
\[
v - u = |\{b^{'}_{h+1}, \dots{} ,b^{'}_{2h}\}\cap \{b_{1}, \dots{} ,b_{2h}\}|\}|.
\]
Define the random variable $Y_{u,v,h,B}(n)$ by
\[
Y_{u,v,h,B}(n) = |\{(b_{1}, \dots{} ,b_{4h-v}): b_{i}\in B, b_{1}, \dots{} ,b_{4h-v}, \textnormal{ distinct } b_{i} < b_{1} \textnormal{ if }
\]
\[
2 \le i \le v \textnormal{ or } i > 2h, b_{1} + \dots{} + b_{2h} = n,
b_{1} + \dots{} + b_{u} + b_{2h+1} + \dots{} + b_{3h-u}
\]
\[
= b_{u+1} + \dots{} + b_{v} + b_{3h-u+1} + \dots{} + b_{4h-v}\}|.
\]
Thus we have
\[
R_{B}^{(2)}(n) \le \sum_{u=1}^{h}\sum_{v=u}^{u+h}Y_{u,v,h,B}(n).
\]
The following lemma gives an upper estimation for $Y_{u,v,h,B}(n)$.

\begin{lem}
With probability 1, for every $1 \le u \le v \le u + h$,
\[
Y_{u,v,h,B}(n) \ll \frac{n^{\frac{1}{4h-1}}}{\log n}.
\]
\end{lem}

It follows from Lemma 7 that with probability 1, there exists a positive integer $C_{u,v,h}(B)$ such that
\[
Y_{u,v,h,B}(n) \le C_{u,v,h}(B)\cdot \frac{n^{\frac{1}{4h-1}}}{\log n}.
\]
Let
\[
C^{(2)}(B) = \sum_{u=1}^{h}\sum_{v=u}^{u+h}C_{u,v,h}(B).
\]
Then with probability 1,
\[
R_{B}^{(2)}(n) \le C^{(2)}(B)\cdot \frac{n^{\frac{1}{4h-1}}}{\log n}.
\]

\noindent \textbf{Case 3.} Upper estimation for $R_{B}^{(3)}(n)$. We have
\[
R_{B}^{(3)}(n) = \sum_{\substack{(d_{1}, \dots{} ,d_{k}, e_{1}, \dots{} ,e_{l})\in (\mathbb{Z}^{+})^{k+l}\\ d_{1} + \dots{} + d_{k} = e_{1} + \dots{} + e_{l} \le h\\ k+l \le 2h - 1}}
|\{(b_{1},\dots{} ,b_{2h}): b_{1},\dots{} ,b_{2h}\in B, b_{1} < \dots{} < b_{2h},
b_{1} + \dots{} + b_{2h} = n,
\]
\[
\exists 1\le i \le 2h,
\exists b_{2}^{'},\dots{} ,b_{k+l}^{'}\in B \textnormal{ distinct }, b_{j}^{'} < b_{i}, d_{1}b_{i} + d_{2}b_{2}^{'} + \dots{} + d_{k}b_{k}^{'} = e_{1}b_{k+1}^{'} + \dots{} + e_{l}b_{k+l}^{'}\}|
\]
\[
= \sum_{\substack{(d_{1}, \dots{} ,d_{k}, e_{1}, \dots{} ,e_{l})\in (\mathbb{Z}^{+})^{k+l}\\ d_{1} + \dots{} + d_{k} = e_{1} + \dots{} + e_{l} \le h\\ k+l \le 2h - 1}}R_{B, d_{1}, \dots{} ,d_{k}, e_{1}, \dots{} ,e_{l}}^{(3)}(n).
\]
For a fixed $(d_{1}, \dots{} ,d_{k}, e_{1}, \dots{} ,e_{l})$, the following lemma gives an upper estimation for the number of such representations.

\begin{lem}
Let $k +l \le 2h - 1$ and $d_{1}, \dots{} ,d_{k}, e_{1}, \dots{} ,e_{l} \in \mathbb{Z}^{+}$ with $d_{1} + \dots{} + d_{k} = e_{1} + \dots{} + e_{l} \le h$.
Then with probability 1, there are only finitely number of solutions of the equation
$d_{1}x_{1} + \dots{} + d_{k}x_{k} = e_{1}x_{k+1} + \dots{} + e_{l}x_{k+l}$, where $x_{1}, \dots{} ,x_{k+l} \in B$ are distinct.
\end{lem}

\noindent According to Lemma 8, almost surely the number of solutions of $d_{1}x_{1} + \dots{} + d_{k}x_{k} = e_{1}x_{k+1} + \dots{} + e_{l}x_{k+l}$, where $x_{1}, \dots{} ,x_{k+l}$ are distinct, in $B$ is at
most $C_{d_{1}, \dots{} ,d_{k}, e_{1}, \dots{} ,e_{l}}^{(3)}(B)$. Then there are at most
\[
\sum_{\substack{(d_{1}, \dots{} ,d_{k}, e_{1}, \dots{} ,e_{l})\in (\mathbb{Z}^{+})^{k+l}\\ d_{1} + \dots{} + d_{k} = e_{1} + \dots{} + e_{l} \le h\\ k+l \le 2h - 1}}C_{d_{1}, \dots{} ,d_{k}, e_{1}, \dots{} ,e_{l}}^{(3)}(B) = C^{(4)}(B)
\]
possibilities for $b_{i}$. By choosing $t = 2h - 1$, $f_{1} = \dots{} = f_{2h-1} = 1$, it follows
from Lemma 6 that the equation $x_{1} + \dots{} + x_{2h-1} = n - b_{i}$, where $x_{1}, \dots{} ,x_{2h-1}$ are distinct, has at most $C^{(5)}(B)$ solutions in the set $B$.
Thus the number of representations in Case 3 is $R_{B}^{(3)}(n) \le C^{(4)}(B)\cdot C^{(5)}(B)$, with probability 1.


We know from (2) that for every $n\ge n_{1}$, $R_{2h,B}(n) \ge C_{h}n^{\frac{1}{4h-1}}$. Then by (4), we have
\[
R_{2h,B\setminus C}(n) \ge C_{h}n^{\frac{1}{4h-1}} - \left(C^{(1)}(B) + C^{(2)}(B)\cdot \frac{n^{\frac{1}{4h-1}}}{\log n} + C^{(4)}(B)\cdot C^{(5)}(B)\right),
\]
which gives that with probability 1, $R_{2h,A}(n) \ge 1$ for every large enough $n$.

Throughout the remaining part of the paper, we prove lemmas 6-8.

\subsection{Proof of Lemma 6}

Let $D$ be an arbitrary set of positive integers. Let
\[
R_{f_{1}, \dots{} ,f_{t}, D}(m) = |\{(d_{1}, \dots{} ,d_{t}): d_{i}\in D, d_{1}, \dots{} ,d_{t} \textnormal{ are distinct, } f_{1}d_{1} + \dots{} + f_{t}d_{t} = m\}|.
\]
Let $R^{\textnormal{Dist}}_{f_{1}, \dots{} ,f_{t}}(m)$ denote the largest nonnegative integer $q$ such that there exist $d^{'}_{1}, \dots{} ,d^{'}_{q_{t}} \in D$ distinct integers such that
$f_{1}d^{'}_{jt+1} + \dots{} + f_{t}d^{'}_{jt+t} = m$ for every $0\le j \le q-1$. Let $E_{f_{1}, \dots{} ,f_{t}}$ denote the event that $R_{f_{1}, \dots{} ,f_{t}, D}(m)$ is not bounded.
Let $E^{Dist}_{f_{1}, \dots{} ,f_{t}}$ denote the event that $R^{\textnormal{Dist}}_{f_{1}, \dots{} ,f_{t}}(m)$ is not bounded. Let
\[
E = \bigcup_{\substack{(f_{1}, \dots{} ,f_{t})\in (\mathbb{Z}^{+})^{n}\\f_{1} + \dots{} + f_{t} \le 2h\\ t\le 2h-1}}E_{f_{1}, \dots{} ,f_{t}}, \hspace*{20mm} E^{\textnormal{Dist}} = \bigcup_{\substack{(f_{1}, \dots{} ,f_{t})\in (\mathbb{Z}^{+})^{n}\\f_{1} + \dots{} + f_{t} \le 2h\\ t\le 2h-1}}E^{\textnormal{Dist}}_{f_{1}, \dots{} ,f_{t}}.
\]
Now we prove that $E = E^{\textnormal{Dist}}$. Obviously, $E^{\textnormal{Dist}}\subseteq E$. Thus we show that if $D$ is an arbitrary subset of the set of positive integers such that
$D\in E$, then $D\in E^{\textnormal{Dist}}$ as well.
We prove that if $2 \le t \le 2h - 1$, $t, u\in \mathbb{Z}^{+}$, then there exists a positive integer $g(t,u)$ such that if there is a positive integer $m$ with $R_{f_{1}, \dots{} ,f_{t}, D}(m) \ge g(t,u)$, then
there exist positive integers $f^{'}_{1}, \dots{} ,f^{'}_{t^{'}}$, $2 \le t^{'} \le 2h - 1$ with $f^{'}_{1} + \dots{} + f^{'}_{t^{'}} \le 2h$ and a positive integer $m^{'}$ such that
$R^{\textnormal{Dist}}_{f^{'}_{1}, \dots{} ,f^{'}_{t}}(m^{'}) \ge u$. We prove by induction on $t$. Assume that $t = 2$. We show that $g(2,u) = 3u - 1$ is suitable. Suppose that $R_{f_{1}, \dots{} ,f_{t}, D}(m)\ge 3u - 1$.
There are at most one representation of the form  $f_{1}d_{1} + f_{2}d_{2} = m$, such that $d_{1} = d_{2}\in D$. It follows that there exist $d^{(1)}_{1}, d^{(1)}_{2}, d^{(2)}_{1}, d^{(2)}_{2}, \dots{} ,d^{(3u-2)}_{1}, d^{(3u-2)}_{2}\in D$ positive integers such that $d^{(i)}_{1}\neq d^{(j)}_{1}$ for every $1 \le i < j \le 3u - 2$ and $f_{1}d^{(i)}_{1} + f_{2}d^{(i)}_{2} = m$ for every $1 \le i \le 3u - 2$. Then it is easy to see that in the representation $f_{1}d^{(i)}_{1} + f_{2}d^{(i)}_{2} = m$ the integers $d^{(i)}_{1}, d^{(i)}_{2}$ appears at most two other representations. In particular, $f_{1}d^{(i)}_{2} + f_{2}\frac{m-f_{1}d^{(i)}_{2}}{f_{2}} = m$ and $f_{1}\frac{m-f_{2}d^{(i)}_{2}}{f_{1}} + f_{2}d^{(i)}_{1} = m$. It follows that there are at least $u$
representations among the $3u - 2$ representations which contain $2u$ distinct integers from $D$.

Assume that the above statement holds for $t - 1 \ge 1$. Now we prove it for $t$. We show that $g(t,u) = t^{2}ug(t - 1,u)$ is suitable. Suppose that $R_{f_{1}, \dots{} ,f_{t}, D}(m)\ge t^{2}ug(t - 1,u) = v$. It follows that there exist $d^{(1)}_{1}, \dots{} ,d^{(1)}_{t}, d^{(2)}_{1}, \dots{} , d^{(2)}_{t}, \dots{} ,d^{(v)}_{1}, \dots{} ,d^{(v)}_{t}\in D$ positive integers such that $d^{(i)}_{1}, \dots{} , d^{(i)}_{t}$ are distinct for $1\le i \le v$
and $f_{1}d^{(i)}_{1} + \dots{} + f_{t}d^{(i)}_{t} = m$ for every $1 \le i \le v$. If there exist a $d\in D$ which appears at least $tg(t - 1,u)$, then there exists a $1\le j \le t$ such that $d^{(i)}_{j} = d$
holds for at least $g(t - 1, u)$ indices $1 \le i \le v$. This implies that $R_{f_{1}, \dots{} ,f_{j-1}, f_{j+1}, \dots{} ,f_{t}, D}(m - f_{j}d) \ge g(t - 1, u)$. It follows from the induction hypothesis that
there exist positive integers $f^{'}_{1}, \dots{} ,f^{'}_{t^{'}}$, $2 \le t^{'} \le 2h - 1$ with $f^{'}_{1} + \dots{} + f^{'}_{t^{'}} \le 2h$ and a positive integer $m^{'}$ such that
$R^{\textnormal{Dist}}_{f^{'}_{1}, \dots{} ,f^{'}_{t^{'}},D}(m^{'}) \ge u$.

If every $d\in D$ appears less than $tg(t - 1,u)$ representations, then the elements from $D$ which appears in the first sum can appear in at most $t^{2}g(t - 1,u)$ sums. Thus among the $v$ representations
there is another one which contains distinct terms from the previous representations.

The elements which appear in these two representations can be at most $2t^{2}g(t - 1,u)$ representations.
Continuing this process, one can get $u$ representations such that any two representations contains distinct terms.
Since every representation contains summands from $D$, we have $R^{Dist}_{f_{1}, \dots{} ,f_{t}, D}(m) \ge u$.

Now we show that if $D\in E$, then $D\in E^{\textnormal{Dist}}$. If $D\in E$, then there exist $f_{1}, \dots{} ,f_{t}$, $2 \le t \le 2h - 1$ positive integers with $f_{1} + \dots{} + f_{t} \le 2h$ such that
$D \in E_{f_{1}, \dots{} ,f_{t}}$. It follows that for every positive integer $u$, there exists a positive integer $m$ such that $R_{f_{1}, \dots{} ,f_{t}, D}(m) \ge g(t,u)$. Thus for every positive integer $u$,
there exist positive integers $f^{'}_{1}, \dots{} ,f^{'}_{t^{'}}$, $2 \le t^{'} \le 2h - 1$ with $f^{'}_{1} + \dots{} + f^{'}_{t} \le 2h$ and a positive integer $m^{'}$ such that $R^{\textnormal{Dist}}_{f^{'}_{1}, \dots{} ,f^{'}_{t},D}(m^{'}) \ge u$. Since there are only finitely many possibilities for $(f^{'}_{1}, \dots{} ,f^{'}_{t^{'}})$, there exist positive integers $f^{'}_{1}, \dots{} ,f^{'}_{t^{'}}$, $2 \le t^{'} \le 2h - 1$ with $f^{'}_{1} + \dots{} + f^{'}_{t} \le 2h$ such that for infinitely many positive integer $u$, there exist a positive integer $m^{'}_{u}$ such that $R^{\textnormal{Dist}}_{f^{'}_{1}, \dots{} ,f^{'}_{t},D}(m^{'}_{u}) \ge u$.
Thus we have $D \in E^{\textnormal{Dist}}_{f^{'}_{1}, \dots{} ,f^{'}_{t}}$ and then $D\in E^{\textnormal{Dist}}$.

To prove Lemma 6 we need to show that $\mathbb{P}(B \in \overline{E_{f_{1}, \dots{} ,f_{t}}}) = 1$, i.e., $\mathbb{P}(B \in E_{f_{1}, \dots{} ,f_{t}}) = 0$. We prove that $\mathbb{P}(B\in E^{\textnormal{Dist}}_{f_{1}, \dots{} ,f_{t}}) = 0$. Let
\[
p_{n} = \mathbb{P}(R^{\textnormal{Dist}}_{f_{1}, \dots{} ,f_{t},D}(n) \ge 4h).
\]
By $t \le 2h - 1$ and (iii) from Lemma 4,
\[
p_{n} \le \left(\sum_{\substack{(x_{1},\dots{} ,x_{t})\in (\mathbb{Z}^{+})^{t}
\\ f_{1}x_{1} + \dots{} + f_{t}x_{t} = n}}
\frac{1}{(x_{1}\cdots{} x_{t})^{\frac{4h-3}{4h-1}}}\right)^{4h}
\]
\[
\ll
\left(\sum_{\substack{(x_{1},\dots{} ,x_{t})\in (\mathbb{Z}^{+})^{t} \\ f_{1}x_{1} + \dots{} + f_{t}x_{t} = n}}
\frac{1}{((f_{1}x_{1})\cdots{} (f_{t}x_{t}))^{\frac{4h-3}{4h-1}}}\right)^{4h}
\]
\[
\ll
\left(\sum_{\substack{(z_{1},\dots{} ,z_{t})\in (\mathbb{Z}^{+})^{t} \\ z_{1} + \dots{} + z_{t} = n}}
\frac{1}{(z_{1}\cdots{} z_{t})^{\frac{4h-3}{4h-1}}}\right)^{4h}
\ll \left(\frac{1}{n^{1-\frac{2t}{4h-1}}}\right)^{4h} \le \frac{1}{n^{\frac{4h}{4h-1}}}.
\]
Since $\frac{4h}{4h-1} > 1$, then by the Borel-Cantelli lemma we get that $R^{\textnormal{Dist}}_{f_{1}, \dots{} ,f_{t},D}(n) \le 4h$ for every large enough $n$, with probability 1.
Thus we have $\mathbb{P}(B\in E^{\textnormal{Dist}}_{f_{1}, \dots{} ,f_{t}}) = 0$ and so $\mathbb{P}(B\in E^{\textnormal{Dist}}) = 0$. Since $E = E^{\textnormal{Dist}}$, we get that
$\mathbb{P}(B\in E) = 0$ and so $\mathbb{P}(B\in E_{f_{1}, \dots{} ,f_{t}}) = 0$ i.e., $\mathbb{P}(B \in \overline{E_{f_{1}, \dots{} ,f_{t}}}) = 1$. The proof of Lemma 6 is completed.

\subsection{Proof of Lemma 8}

Let $d_{1}, \dots{} ,d_{k}, e_{1}, \dots{} ,e_{l}$, $k + l \le 2h - 1$ be positive integers with $d_{1} + \dots{} + d_{k} =  e_{1} + \dots{}  + e_{l} \le h$. Let $n$ be a positive integer.
Let
\[
q_{n} = \mathbb{P}(\exists b_{1}, \dots{} ,b_{k+l}\in B \textnormal{ distinct }, d_{1}b_{1} + \dots{} + d_{k}b_{k} = e_{1}b_{k+1} + \dots{} + e_{l}b_{k+l} = n).
\]
By (iii) in Lemma 4,
\[
q_{n} \le \sum_{\substack{(x_{1},\dots{} ,x_{k+l})\in (\mathbb{Z}^{+})^{k+l} \\ \textnormal{the $x_{i}$'s are distinct} \\ d_{1}x_{1} + \dots{} + d_{k}x_{k} = e_{1}x_{k+1} + \dots{} + e_{l}x_{k+l} = n}}
\frac{1}{(x_{1}\cdots{} x_{k+l})^{\frac{4h-3}{4h-1}}}
\]
\[
\le \left(\sum_{\substack{(x_{1},\dots{} ,x_{k})\in (\mathbb{Z}^{+})^{k} \\ \textnormal{the $x_{i}$'s are distinct} \\ d_{1}x_{1} + \dots{} + d_{k}x_{k} = n}}
\frac{1}{(x_{1}\cdots{} x_{k})^{\frac{4h-3}{4h-1}}}\right)\cdot \left(\sum_{\substack{(x_{k+1},\dots{} ,x_{k+l})\in (\mathbb{Z}^{+})^{l} \\ \textnormal{the $x_{i}$'s are distinct} \\ e_{1}x_{k+1} + \dots{} + e_{l}x_{k+l} = n}}
\frac{1}{(x_{k+1}\cdots{} x_{k+l})^{\frac{4h-3}{4h-1}}}\right)
\]
\[
\ll \left(\sum_{\substack{(x_{1},\dots{} ,x_{k})\in (\mathbb{Z}^{+})^{k} \\ \textnormal{the $x_{i}$'s are distinct} \\ d_{1}x_{1} + \dots{} + d_{k}x_{k} = n}}
\frac{1}{((d_{1}x_{1})\cdots{} (d_{k}x_{k}))^{\frac{4h-3}{4h-1}}}\right)\cdot \left(\sum_{\substack{(x_{k+1},\dots{} ,x_{k+l})\in (\mathbb{Z}^{+})^{l} \\ \textnormal{the $x_{i}$'s are distinct} \\ e_{1}x_{k+1} + \dots{} + e_{l}x_{k+l} = n}}
\frac{1}{((e_{1}x_{k+1})\cdots{} (e_{l}x_{k+l}))^{\frac{4h-3}{4h-1}}}\right)
\]
\[
\le \left(\sum_{\substack{(z_{1},\dots{} ,z_{k})\in (\mathbb{Z}^{+})^{k} \\ z_{1} + \dots{} + z_{k} = n}}
\frac{1}{(z_{1}\cdots{} z_{k})^{\frac{4h-3}{4h-1}}}\right)\cdot \left(\sum_{\substack{(z_{k+1},\dots{} ,z_{k+l})\in (\mathbb{Z}^{+})^{l} \\ z_{k+1} + \dots{} + z_{k+l} = n}}
\frac{1}{(z_{k+1}\cdots{} z_{k+l})^{\frac{4h-3}{4h-1}}}\right)
\]
\[
\ll \frac{1}{n^{1-\frac{2k}{4h-1}}}\cdot \frac{1}{n^{1-\frac{2l}{4h-1}}} = \frac{1}{n^{2-\frac{2(k+l)}{4h-1}}}.
\]
Since $k + l \le 2h - 1$, then $2 - \frac{2(k+l)}{4h-1} > 1$. It follows that the infinite series $\sum_{n=1}^{\infty}q_{n}$ is convergent, thus by the Borel-Cantelli lemma almost surely there exist only
finitely many positive integer $n$ for which there exist distinct positive integers $b_{1}, \dots{} ,b_{k+l} \in B$ such that $d_{1}b_{1} + \dots{} + d_{k}b_{k} = e_{1}b_{k+1} + \dots{} + e_{l}b_{k+l}$.
The proof of Lemma 8 is completed.

\subsection{Proof of Lemma 7}

Let $\overline{\alpha} = (\alpha_{1}, \dots{} ,\alpha_{n})$, $\alpha_{i}\in \mathbb{N}$, $\alpha_{i} \ge 0$. We will prove Lemma 7 by using Lemma 2. Recall
\[
t_{n} =
\left\{
\begin{aligned}
1 \textnormal{, if } n \in \mathcal{A} \\
0 \textnormal{, if } n \notin \mathcal{A}
\end{aligned} \hspace*{3mm}.
\right.
\]
The random variable $Y_{u,v,h,B}(n)$ is clearly
\begin{align*}
Y_{u,v,h,B}(n) = &\left(\sum_{\substack{(x_{1},\dots{} ,x_{u})\in (\mathbb{Z}^{+})^{u} \\ \textnormal{ $x_{1}, \dots{} ,x_{u}$ are distinct} \\ x_{i}< x_{1} \textnormal{, if } 2\le i \le u}}t_{x_{1}}\cdots{} t_{x_{u}}\left(\sum_{\substack{(x_{u+1},\dots{} ,x_{v})\in (\mathbb{Z}^{+})^{v-u} \\ \textnormal{ $x_{1}, \dots{} ,x_{v}$ are distinct} \\ x_{i}< x_{1} \textnormal{, if } u < i \le v}}t_{x_{u+1}}\cdots{} t_{x_{v}}
\right.\right.\\
&\mkern10mu\left(\sum_{\substack{(x_{v+1},\dots{} ,x_{2h})\in (\mathbb{Z}^{+})^{2h-v} \\ \textnormal{ $x_{1}, \dots{} ,x_{2h}$ are distinct} \\ x_{1} + \dots{} + x_{2h} = n}}t_{x_{v+1}}\cdots{} t_{x_{2h}}\left(\sum_{\substack{(x_{2h+1},\dots{} ,x_{3h-u})\in (\mathbb{Z}^{+})^{h-u}\\ \textnormal{ $x_{1}, \dots{} ,x_{3h-u}$ are distinct} \\ x_{i}< x_{1} \textnormal{, if } 2h+1\le i \le 3h-u}}t_{x_{2h+1}}\cdots{} t_{x_{3h-u}}\right.\right.\\
&\left.\left.\left.\left.\left(\sum_{\substack{(x_{3h-u+1},\dots{} ,x_{4h-v})\in (\mathbb{Z}^{+})^{h-(v-u)}\\ \textnormal{ $x_{1}, \dots{} ,x_{4h-v}$ are distinct} \\ x_{i}< x_{1} \textnormal{, if } 3h-u+1\le i \le 4h-v}}t_{x_{3h-u+1}}\cdots{} t_{x_{4h-v}}\right)\right)\right)\right)\right).
\end{align*}
Let us denote $\partial^{(\overline{\alpha})}Y_{u,v,h,B}(n) = Y^{(\overline{\alpha})}_{u,v,h,B}(n)$. We show that
\begin{equation}
\mathbb{E}(Y^{(\overline{\alpha})}_{u,v,h,B}(n)) = O\left(\frac{n^{\frac{1}{4h-1}}}{(\log n)^{4h}}\right),
\end{equation}
where the constant depends only on $h$. It follows from (5) that
\[
\mathbb{E}(Y_{u,v,h,B}(n)) = O\left(\frac{n^{\frac{1}{4h-1}}}{(\log n)^{4h}}\right) \hspace*{10mm} \mathbb{E}_{\ge 0}(Y^{(\overline{\alpha})}_{u,v,h,B}(n)) = O\left(\frac{n^{\frac{1}{4h-1}}}{(\log n)^{4h}}\right)
\]
\[
\mathbb{E}_{\ge 1}(Y^{(\overline{\alpha})}_{u,v,h,B}(n)) = O\left(\frac{n^{\frac{1}{4h-1}}}{(\log n)^{4h}}\right).
\]
Applying Lemma 2  with $k = 4h - v \le 4h - 1$ and $\lambda = 16h\log n$, we get that
\[
\mathbb{P}\left(|Y_{u,v,h,B}(n) - \mathbb{E}(Y_{u,v,h,B}(n))| \ge C_{4h-v}(16h\log n)^{4h-v-\frac{1}{2}}\sqrt{\mathbb{E}_{\ge 0}(Y^{(\overline{\alpha})}_{u,v,h,B}(n))\cdot \mathbb{E}_{\ge 1}(Y^{(\overline{\alpha})}_{u,v,h,B}(n))}\right)
\]
\[
= O\left(e^{-\frac{16h\log n}{4}+(4h-2)\log n}\right) = O\left(\frac{1}{n^{2}}\right).
\]

Then by the Borel-Cantelli lemma, we get that with probability 1,

\[
|Y_{u,v,h,B}(n) - \mathbb{E}(Y_{u,v,h,B}(n))| < C_{4h-v}(16h\log n)^{4h-v-\frac{1}{2}}\sqrt{\mathbb{E}_{\ge 0}(Y^{(\overline{\alpha})}_{u,v,h,B}(n))\cdot \mathbb{E}_{\ge 1}(Y^{(\overline{\alpha})}_{u,v,h,B}(n))}
\]
holds for every $n$ large enough. Then with probability 1,
\[
Y_{u,v,h,B}(n) = O\left((16h\log n)^{4h-v-\frac{1}{2}}\cdot \frac{n^{\frac{1}{4h-1}}}{(\log n)^{4h}}\right) = O\left(\frac{n^{\frac{1}{4h-1}}}{\log n}\right).
\]
Thus it is enough to prove (5). Obviously,
\begin{align*}
Y_{u,v,h,B}^{(\overline{\alpha})}(n) = \frac{\partial}{\partial \overline{\alpha}}&\left(\sum_{\substack{(x_{1},\dots{} ,x_{u})\in (\mathbb{Z}^{+})^{u} \\ \textnormal{ $x_{1}, \dots{} ,x_{u}$ are distinct} \\ x_{i}< x_{1} \textnormal{, if } 2\le i \le u}}t_{x_{1}}\cdots{} t_{x_{u}}\left(\sum_{\substack{(x_{u+1},\dots{} ,x_{v})\in (\mathbb{Z}^{+})^{v-u} \\ \textnormal{ $x_{1}, \dots{} ,x_{v}$ are distinct} \\ x_{i}< x_{1} \textnormal{, if } u < i \le v}}t_{x_{u+1}}\cdots{} t_{x_{v}}
\right.\right.\\
&\mkern10mu\left(\sum_{\substack{(x_{v+1},\dots{} ,x_{2h})\in (\mathbb{Z}^{+})^{2h-v} \\ \textnormal{ $x_{1}, \dots{} ,x_{2h}$ are distinct} \\ x_{1} + \dots{} + x_{2h} = n}}t_{x_{v+1}}\cdots{} t_{x_{2h}}\left(\sum_{\substack{(x_{2h+1},\dots{} ,x_{3h-u})\in (\mathbb{Z}^{+})^{h-u}\\ \textnormal{ $x_{1}, \dots{} ,x_{3h-u}$ are distinct} \\ x_{i}< x_{1} \textnormal{, if } 2h+1\le i \le 3h-u}}t_{x_{2h+1}}\cdots{} t_{x_{3h-u}}\right.\right.\\
&\left.\left.\left.\left.\left(\sum_{\substack{(x_{3h-u+1},\dots{} ,x_{4h-v})\in (\mathbb{Z}^{+})^{h-(v-u)}\\ \textnormal{ $x_{1}, \dots{} ,x_{4h-v}$ are distinct} \\ x_{i}< x_{1} \textnormal{, if } 3h-u+1\le i \le 4h-v}}t_{x_{3h-u+1}}\cdots{} t_{x_{4h-v}}\right)\right)\right)\right)\right).
\end{align*}
We can assume that $\alpha_{i}\in \{0,1\}$, otherwise $Y_{u,v,h,B}^{(\overline{\alpha})}(n) = 0$ because the integers $x_{1}, \dots{} ,x_{4h-v}$ are distinct. Let $|\overline{\alpha}| = \sum_{i=1}^{n}\alpha_{i}$,
where we can assume that $|\overline{\alpha}|\le 4h-v$, otherwise $Y_{u,v,h,B}^{(\overline{\alpha})}(n) = 0$.
Let $\overline{\alpha} = \overline{\beta} + \overline{\gamma} + \overline{\delta} + \overline{\mu} + \overline{\nu}$, where $\overline{\beta} = (\beta_{1}, \dots{} ,\beta_{n})$, $\overline{\gamma} = (\gamma_{1}, \dots{} ,\gamma_{n})$,
$\overline{\delta} = (\delta_{1}, \dots{} ,\delta_{n})$, $\overline{\mu} = (\mu_{1}, \dots{} ,\mu_{n})$, $\overline{\nu} = (\nu_{1}, \dots{} ,\nu_{n})$ and $\beta_{i}, \gamma_{i}, \delta_{i}, \mu_{i}, \nu_{i}\in \{0,1\}$ for $1\le i \le n$.
Then we have
\par
\bigskip
$
\displaystyle Y^{\overline{(\alpha)}}_{u,v,h,B}(n) =
$
\begin{align*}
&= \sum_{\substack{(\overline{\beta}, \overline{\gamma}, \overline{\delta}, \overline{\mu}, \overline{\nu})\\ \overline{\alpha} = \overline{\beta} + \overline{\gamma} + \overline{\delta} + \overline{\mu} + \overline{\nu}}}\left(\sum_{\substack{(x_{1},\dots{} ,x_{u})\in (\mathbb{Z}^{+})^{u} \\ \textnormal{ $x_{1}, \dots{} ,x_{u}$ are distinct} \\ x_{i}< x_{1} \textnormal{, if } 2\le i \le u}}\left(\frac{\partial}{\partial \overline{\beta}}t_{x_{1}}\cdots{} t_{x_{u}}\right)\left(\sum_{\substack{(x_{u+1},\dots{} ,x_{v})\in (\mathbb{Z}^{+})^{v-u} \\ \textnormal{the $x_{1}, \dots{} ,x_{v}$ are distinct} \\ x_{i}< x_{1} \textnormal{, if } u < i \le v}}\left(\frac{\partial}{\partial \overline{\gamma}}t_{x_{u+1}}\cdots{} t_{x_{v}}\right)\right.\right.\\
&\mkern10mu\left(\sum_{\substack{(x_{v+1},\dots{} ,x_{2h})\in (\mathbb{Z}^{+})^{2h-v} \\ \textnormal{$x_{1}, \dots{} ,x_{2h}$ are distinct} \\ x_{1} + \dots{} + x_{2h} = n}}\left(\frac{\partial}{\partial \overline{\delta}}t_{x_{v+1}}\cdots{} t_{x_{2h}}\right)\left(\sum_{\substack{(x_{2h+1},\dots{} ,x_{3h-u})\in (\mathbb{Z}^{+})^{h-u}\\ \textnormal{ $x_{1}, \dots{} ,x_{3h-u}$ are distinct} \\ x_{i}< x_{1} \textnormal{, if } 2h+1\le i \le 3h-u}}\left(\frac{\partial}{\partial \overline{\mu}}t_{x_{2h+1}}\cdots{} t_{x_{3h-u}}\right)
\right.\right.\\
&\mkern10mu\left.\left.\left.\left.\left(\sum_{\substack{(x_{3h-u+1},\dots{} ,x_{4h-v})\in (\mathbb{Z}^{+})^{h+u-v}\\ \textnormal{ $x_{1}, \dots{} ,x_{4h-v}$ are distinct} \\ x_{i}< x_{1} \textnormal{, if } 3h-u+1\le i \le 4h-v}}\left(\frac{\partial}{\partial \overline{\nu}}t_{x_{3h-u+1}}\cdots{} t_{x_{4h-v}}\right)\right)\right)\right)\right)\right)\\
&= \sum_{(\overline{\beta}, \overline{\gamma}, \overline{\delta}, \overline{\mu}, \overline{\nu})}Y^{(\overline{\beta}, \overline{\gamma}, \overline{\delta}, \overline{\mu}, \overline{\nu})}_{u,v,h,B}(n).\\
\end{align*}
It is clear that there are at most $5^{|\overline{\alpha}|}\le 5^{4h-v}$ possibilities for the $5$-tuples $(\overline{\beta}, \overline{\gamma}, \overline{\delta}, \overline{\mu}, \overline{\nu})$.
It is enough to prove that for a fixed $5$-tuples
$(\overline{\beta}, \overline{\gamma}, \overline{\delta}, \overline{\mu}, \overline{\nu})$,
\[
\mathbb{E}(Y^{(\overline{\beta}, \overline{\gamma}, \overline{\delta}, \overline{\mu}, \overline{\nu})}_{u,v,h,B}(n)) = O\left(\frac{n^{\frac{1}{4h-1}}}{(\log n)^{4h}}\right).
\]
Assume that the nonzero $\beta_{i}, \gamma_{i}, \delta_{i}, \mu_{i}, \nu_{i}$'s have the indices $i_{1}, \dots{} ,i_{|\overline{\beta}|}$ and $i_{|\overline{\beta}|+1}, \dots{} ,i_{|\overline{\beta}|+|\overline{\gamma}|}$ and
$i_{|\overline{\beta}|+|\overline{\gamma}|+1}, \dots{} ,i_{|\overline{\beta}|+|\overline{\gamma}|+|\overline{\delta}|}$ and $i_{|\overline{\beta}|+|\overline{\gamma}|+|\overline{\delta}|+1}, \dots{} ,i_{|\overline{\beta}|+|\overline{\gamma}|+|\overline{\delta}|+|\overline{\mu}|}$ and
$i_{|\overline{\beta}|+|\overline{\gamma}|+|\overline{\delta}|+|\overline{\mu}|+1}, \dots{}$

,$i_{|\overline{\beta}|+|\overline{\gamma}|+|\overline{\delta}|+|\overline{\mu}|+|\overline{\nu}|}$, respectively.
If $i_{1} + \dots{} + i_{|\overline{\beta}|+|\overline{\gamma}|+|\overline{\delta}|} > n$, then $Y^{(\overline{\beta}, \overline{\gamma}, \overline{\delta}, \overline{\mu}, \overline{\nu})}_{u,v,h,B}(n) = 0$. If $i_{1} + \dots{} + i_{|\overline{\beta}|+|\overline{\gamma}|+|\overline{\delta}|} = n$, then we can assume that $|\overline{\beta}| = u$, $|\overline{\gamma}| = v - u$ and $|\overline{\delta}| = 2h - v$, otherwise $Y^{(\overline{\beta}, \overline{\gamma}, \overline{\delta}, \overline{\mu}, \overline{\nu})}_{u,v,h,B}(n) = 0$. Then we have
\[
Y^{(\overline{\beta}, \overline{\gamma}, \overline{\delta}, \overline{\mu}, \overline{\nu})}_{u,v,h,B}(n) \le |\{(w_{1}, \dots{} ,w_{2h-v}): w_{i}\in B, 1\le w_{i} \le n, w_{1}, \dots{} ,w_{2h-v} \textnormal{ are distinct, }
\]
\[
w_{1} + \dots{} + w_{h-u} - (w_{h-u+1} + \dots{} + w_{2h-v}) =
i_{u+1} + \dots{} + i_{v} - (i_{1} + \dots{} + i_{u}) = M\}|.
\]
It follows from (iv) in Lemma 4 that
\[
\mathbb{E}(Y^{(\overline{\beta}, \overline{\gamma}, \overline{\delta}, \overline{\mu}, \overline{\nu})}_{u,v,h,B}(n))\le \sum_{\substack{(w_{1},\dots{} ,w_{2h-v})\in (\mathbb{Z}^{+})^{2h-v} \\ w_{1} + \dots{} + w_{h-u} - (w_{h-u+1} + \dots{} + w_{2h-v}) =
i_{u+1} + \dots{} + i_{v} - (i_{1} + \dots{} + i_{u}) = M}}\frac{1}{(w_{1} \cdots{} w_{2h-v})^{\frac{4h-3}{4h-1}}}
\]
\[
\ll \frac{1}{(|M|+1)^{1-\frac{2(2h-v)}{4h-1}}} = \frac{1}{(|M|+1)^{\frac{2v-1}{4h-1}}} = O\left(\frac{n^{\frac{1}{4h-1}}}{(\log n)^{4h}}\right).
\]
Now, we assume that $i_{1} + \dots{} + i_{|\overline{\beta}|+|\overline{\gamma}|+|\overline{\delta}|} < n$. We denote the following conditions by $\textnormal{cond}_{3}$ and $\textnormal{cond}_{4}$, respectively:
\[
\textnormal{cond}_{3}: \hspace*{20mm} z_{1} + \dots{} + z_{2h-|\overline{\beta}|-|\overline{\gamma}|-|\overline{\delta}|} = n - (i_{1} + \dots{} + i_{|\overline{\beta}|+|\overline{\gamma}|+|\overline{\delta}|}),
\]
\[
\textnormal{cond}_{4}: \hspace*{5mm} z_{1} + \dots{} + z_{u-|\overline{\beta}|} + z_{2h-|\overline{\beta}|-|\overline{\gamma}|-|\overline{\delta}|+1} + \dots{} + z_{3h-u-|\overline{\beta}|-|\overline{\gamma}|-|\overline{\delta}|-|\overline{\mu}|}
\]
\[
- (z_{u-|\overline{\beta}|+1} + \dots{} + z_{v-|\overline{\beta}|-|\overline{\gamma}|} + z_{3h-u-|\overline{\beta}|-|\overline{\gamma}|-|\overline{\delta}|-|\overline{\mu}|+1}+ \dots{} + z_{4h-v-|\overline{\beta}|-|\overline{\gamma}|-|\overline{\delta}|-|\overline{\mu}|-|\overline{\nu}|})
\]
\[
= i_{|\overline{\beta}|+1} + \dots{} + i_{|\overline{\beta}|+|\overline{\gamma}|} + i_{|\overline{\beta}|+|\overline{\gamma}|+|\overline{\delta}|+|\overline{\mu}|+1} + \dots{}
+ i_{|\overline{\beta}|+|\overline{\gamma}|+|\overline{\delta}|+|\overline{\mu}|+|\overline{\nu}|}
\]
\[
- (i_{1} + \dots{} + i_{|\overline{\beta}|} + i_{|\overline{\beta}|+|\overline{\gamma}|+|\overline{\delta}|+1} + \dots{} + i_{|\overline{\beta}|+|\overline{\gamma}|+|\overline{\delta}|+|\overline{\mu}|}) = M.
\]
Then we have
\begin{align*}
Y^{(\overline{\beta}, \overline{\gamma}, \overline{\delta}, \overline{\mu}, \overline{\nu})}_{u,v,h,B}(n) &= |\{(z_{1}, \dots{} ,z_{4h-v-|\overline{\beta}|-|\overline{\gamma}|-|\overline{\delta}|-|\overline{\mu}|-|\overline{\nu}|}): z_{i}\in B, z_{j} \neq i_{t}\\
&\textnormal{for} \hspace*{2mm}
1 \le j \le 4h-v-|\overline{\beta}|-|\overline{\gamma}|-|\overline{\delta}|-|\overline{\mu}|-|\overline{\nu}|,
1 \le t \le |\overline{\beta}|+|\overline{\gamma}|+|\overline{\delta}|+|\overline{\mu}|+|\overline{\nu}|,\\
&\textnormal{$cond_{3}$ and $cond_{4}$ hold}\}|.\\
\end{align*}
It follows that
\begin{equation}
\mathbb{E}(Y^{(\overline{\beta}, \overline{\gamma}, \overline{\delta}, \overline{\mu}, \overline{\nu})}_{u,v,h,B}(n)) \le
\end{equation}
\[
\le \sum_{\substack{(z_{1}, \dots{} ,z_{2h-|\overline{\beta}|-|\overline{\gamma}|-|\overline{\delta}|)}\\ \textnormal{$cond_{3}$ holds}}}\left(\frac{1}{(z_{1} \cdots{} z_{2h-|\overline{\beta}|-|\overline{\gamma}|-|\overline{\delta}|})^{\frac{4h-3}{4h-1}}}\cdot\right.
\]
\[
\left.\cdot\sum_{\substack{(z_{2h-|\overline{\beta}|-|\overline{\gamma}|-|\overline{\delta}|+1},\dots{}, z_{4h-v-|\overline{\beta}|-|\overline{\gamma}|-|\overline{\delta}|-|\overline{\mu}|-|\overline{\nu}|})\\ \textnormal{$cond_{4}$ holds}}}
\frac{1}{(z_{2h-|\overline{\beta}|-|\overline{\gamma}|-|\overline{\delta}|+1} \cdots{} z_{4h-v-|\overline{\beta}|-|\overline{\gamma}|-|\overline{\delta}|-|\overline{\mu}|-|\overline{\nu}|})^{\frac{4h-3}{4h-1}}}\right).
\]
By (iv) of Lemma 4, the inner sum
\[
\sum_{\substack{(z_{2h-|\overline{\beta}|-|\overline{\gamma}|-|\overline{\delta}|+1},\dots{}, z_{4h-v-|\overline{\beta}|-|\overline{\gamma}|-|\overline{\delta}|-|\overline{\mu}|-|\overline{\nu}|})\\ \textnormal{$cond_{4}$ holds}}}\frac{1}{(z_{2h-|\overline{\beta}|-|\overline{\gamma}|-|\overline{\delta}|+1} \cdots{} z_{4h-v-|\overline{\beta}|-|\overline{\gamma}|-|\overline{\delta}|-|\overline{\mu}|-|\overline{\nu}|})^{\frac{4h-3}{4h-1}}}
\]
\[
\ll \frac{1}{(1+|M|)^{1-\frac{2(2h-v-|\overline{\mu}|-|\overline{\nu}|)}{4h-1}}}
\ll \frac{1}{(1+|M|)^{\frac{2v-1}{4h-1}}} = O(1).
\]
On the other hand, by (iii) of Lemma 4,
\[
\sum_{\substack{(z_{1}, \dots{} ,z_{2h-|\overline{\beta}|-|\overline{\gamma}|-|\overline{\delta}|}) \\ \textnormal{$cond_{3}$ holds}}}\frac{1}{(z_{1} \cdots{} z_{2h-|\overline{\beta}|-|\overline{\gamma}|-|\overline{\delta}|})^{\frac{4h-3}{4h-1}}}
\ll \frac{1}{(n - (i_{1} + \dots{} + i_{|\overline{\beta}|+|\overline{\gamma}|+|\overline{\delta}|}))^{1-\frac{2(2h-|\overline{\beta}|-|\overline{\gamma}|-|\overline{\delta}|)}{4h-1}}}
\]
\[
= \frac{1}{(n - (i_{1} + \dots{} + i_{|\overline{\beta}|+|\overline{\gamma}|+|\overline{\delta}|}))^{\frac{2(|\overline{\beta}|+|\overline{\gamma}|+|\overline{\delta}|)-1}{4h-1}}}.
\]
Thus if $|\overline{\beta}|+|\overline{\gamma}|+|\overline{\delta}| \ge 1$, then  $\mathbb{E}(Y^{(\overline{\beta}, \overline{\gamma}, \overline{\delta}, \overline{\mu}, \overline{\nu})}_{u,v,h,B}(n)) = O(1)$.
Now we can assume that
$|\overline{\beta}| = |\overline{\gamma}| = |\overline{\delta}| = 0$, i.e., $\overline{\beta} = \overline{\gamma} = \overline{\delta} = \overline{0}$. Then the equation in $cond_{4}$ can be written in the form
\[
z_{1} + \dots{} + z_{u} + z_{2h+1} + \dots{} + z_{3h-u-|\overline{\mu}|}
- (z_{u+1} + \dots{} + z_{v} + z_{3h-u-|\overline{\mu}|+1}+ \dots{} + z_{4h-v-|\overline{\mu}|-|\overline{\nu}|})
\]
\[
= i_{|\overline{\mu}|+1} + \dots{} + i_{|\overline{\mu}|+|\overline{\nu}|} - (i_{1} + \dots{} + i_{|\overline{\mu}|}),
\]
that is, $cond_{4}$ is equivalent to
\[
\textnormal{$cond_{5}$}: \hspace*{5mm} z_{2h+1} + \dots{} + z_{3h-u-|\overline{\mu}|}
- (z_{3h-u-|\overline{\mu}|+1}+ \dots{} + z_{4h-v-|\overline{\mu}|-|\overline{\nu}|}) = z_{u+1} + \dots{} + z_{v} - (z_{1} + \dots{} + z_{u})
\]
\[
+ i_{|\overline{\mu}|+1} + \dots{} + i_{|\overline{\mu}|+|\overline{\nu}|} - (i_{1} + \dots{} + i_{|\overline{\mu}|}).
\]
If $K = i_{1} + \dots{} + i_{|\overline{\mu}|} - (i_{|\overline{\mu}|+1} + \dots{} + i_{|\overline{\mu}|+|\overline{\nu}|})$, then, by replacing $cond_{4}$ to $cond_{5}$ and by (iv) of Lemma 4,
the inner sum in (6) can be estimated in the following way.
\[
\sum_{\substack{(z_{2h+1},\dots{}, z_{4h-v-|\overline{\mu}|-|\overline{\nu}|})\\ \textnormal{$cond_{5}$ holds}}}\frac{1}{(z_{2h+1} \cdots{} z_{4h-v-|\overline{\mu}|-|\overline{\nu}|})^{\frac{4h-3}{4h-1}}} \ll
\]
\[
\ll \frac{1}{(|z_{1} + \dots{} + z_{u} - (z_{u+1} + \dots{} + z_{v}) + K| + 1)^{\frac{2v-1}{4h-1}}},
\]
and then
\[
\mathbb{E}(Y^{(\overline{\beta}, \overline{\gamma}, \overline{\delta}, \overline{\mu}, \overline{\nu})}_{u,v,h,B})(n) \ll
\]
\[
\ll \sum_{\substack{(z_{1},\dots{} ,z_{2h})\in (\mathbb{Z}^{+})^{2h} \\ z_{1} + \dots{} + z_{2h} = n}}\frac{1}{(z_{1} \cdots{} z_{2h})^{\frac{4h-3}{4h-1}}}\cdot \frac{1}{(|z_{1} + \dots{} + z_{u} - (z_{u+1} + \dots{} + z_{v}) + K| + 1)^{\frac{2v-1}{4h-1}}}.
\]
We will prove that
\[
\sum_{\substack{(z_{1},\dots{} ,z_{2h})\in (\mathbb{Z}^{+})^{2h} \\ z_{1} + \dots{} + z_{2h} = n}}\frac{1}{(z_{1} \cdots{} z_{2h})^{\frac{4h-3}{4h-1}}}\cdot \frac{1}{(|z_{1} + \dots{} + z_{u} - (z_{u+1} + \dots{} + z_{v}) + K| + 1)^{\frac{2v-1}{4h-1}}}
\]
\[
= O\left(\frac{n^{\frac{1}{4h-1}}}{(\log n)^{4h}}\right).
\]
We have three cases.

\textbf{Case 1.} $v = 2h$. Then $u = h$ and $\overline{\mu} = \overline{\nu} = \overline{0}$, and so
\[
Y^{(\overline{\beta}, \overline{\gamma}, \overline{\delta}, \overline{\mu}, \overline{\nu})}_{u,v,h,B}(n) = |\{(z_{1}, \dots{} ,z_{2h}): z_{i}\in B, z_{1} + \dots{} + z_{2h} = n, z_{1} + \dots{} + z_{h} = z_{h+1} + \dots{}
+ z_{2h}, z_{i}^{'}s
\]
\[
\textnormal{are distinct}, z_{i} < z_{1}, \textnormal{ if } 2 \le i \le 2h\}|.
\]
It follows from (iii) of Lemma 4 that
\[
\mathbb{E}(Y^{(\overline{\beta}, \overline{\gamma}, \overline{\delta}, \overline{\mu}, \overline{\nu})}_{u,v,h,B})(n) \ll \left(\sum_{\substack{(z_{1},\dots{} ,z_{h})\in (\mathbb{Z}^{+})^{h} \\ z_{1} + \dots{} + z_{h} = \frac{n}{2}}}\frac{1}{(z_{1} \cdots{} z_{h})^{\frac{4h-3}{4h-1}}}\right)^{2} \ll \left(\frac{1}{n^{1-\frac{2h}{4h-1}}}\right)^{2} = O(1).
\]

\textbf{Case 2.} $u = v$. Then $u = v \le h$ and so
\bigskip
\par
$\mathbb{E}(Y^{(\overline{\beta}, \overline{\gamma}, \overline{\delta}, \overline{\mu}, \overline{\nu})}_{u,v,h,B})(n) \ll$

\begin{align*}
&\ll \sum_{\substack{(z_{1},\dots{} ,z_{2h})\in (\mathbb{Z}^{+})^{2h} \\ z_{1} + \dots{} + z_{2h} = n}}\frac{1}{(z_{1} \cdots{} z_{2h})^{\frac{4h-3}{4h-1}}}\cdot \frac{1}{(|z_{1} + \dots{} + z_{u} - (z_{u+1} + \dots{} + z_{v}) + K| + 1)^{\frac{2v-1}{4h-1}}}\\
&= \sum_{m=1}^{n-1}\left(\sum_{\substack{(z_{1},\dots{} ,z_{u}) \\ z_{1} + \dots{} + z_{u} = m}}\left(\frac{1}{(z_{1} \cdots{} z_{u})^{\frac{4h-3}{4h-1}}}\cdot \frac{1}{(|z_{1} + \dots{} + z_{u} + K| + 1)^{\frac{2v-1}{4h-1}}}\cdot\right.\right.\\
&\left.\left.\cdot\sum_{\substack{(z_{u+1},\dots{} ,z_{2h}) \\ z_{u+1} + \dots{} + z_{2h} = n-m}}\frac{1}{(z_{u+1} \cdots{} z_{2h})^{\frac{4h-3}{4h-1}}}\right)\right).\\
\end{align*}
It follows from (iii) of Lemma 4 that
\[
\mathbb{E}(Y^{(\overline{\beta}, \overline{\gamma}, \overline{\delta}, \overline{\mu}, \overline{\nu})}_{u,v,h,B})(n)
\ll \sum_{m=1}^{n-1}\left(\left(\sum_{\substack{(z_{1},\dots{} ,z_{u}) \\ z_{1} + \dots{} + z_{u} = m}}\frac{1}{(z_{1} \cdots{} z_{u})^{\frac{4h-3}{4h-1}}}\right)\cdot \frac{1}{(|m + K| + 1)^{\frac{2u-1}{4h-1}}}\cdot \frac{1}{(n - m)^{\frac{2u-1}{4h-1}}}
\right)
\]
\[
\ll \sum_{m=1}^{n-1}\frac{1}{m^{1-\frac{2u}{4h-1}}} \cdot \frac{1}{(|m + K| + 1)^{\frac{2u-1}{4h-1}}}\cdot \frac{1}{(n - m)^{\frac{2u-1}{4h-1}}}.
\]
We have three subcases.

\textbf{Subase 2.1.} $K \ge 0$. It follows from (i) of Lemma 4 that

\[
\mathbb{E}(Y^{(\overline{\beta}, \overline{\gamma}, \overline{\delta}, \overline{\mu}, \overline{\nu})}_{u,v,h,B}(n)) \ll \sum_{m=1}^{n-1}\frac{1}{m^{1-\frac{2u}{4h-1}}} \cdot \frac{1}{m^{\frac{2u-1}{4h-1}}}\cdot \frac{1}{(n - m)^{\frac{2u-1}{4h-1}}} = \sum_{m=1}^{n-1}\frac{1}{m^{1-\frac{1}{4h-1}}} \cdot \frac{1}{(n - m)^{\frac{2u-1}{4h-1}}}
\]
\[
\ll \frac{1}{m^{\frac{2u-2}{4h-1}}} = O(1).
\]

\textbf{Subase 2.2.} $-\frac{n}{3} < K < 0$. Then we have

\[
\mathbb{E}(Y^{(\overline{\beta}, \overline{\gamma}, \overline{\delta}, \overline{\mu}, \overline{\nu})}_{u,v,h,B}(n)) \ll \sum_{m=1}^{-\lfloor K/2 \rfloor}\frac{1}{m^{1-\frac{2u}{4h-1}}} \cdot \frac{1}{(|m + K| + 1)^{\frac{2u-1}{4h-1}}}\cdot \frac{1}{(n - m)^{\frac{2u-1}{4h-1}}}
\]
\[
+ \sum_{m=-\lfloor K/2 \rfloor+1}^{-K}\frac{1}{m^{1-\frac{2u}{4h-1}}} \cdot \frac{1}{(|m + K| + 1)^{\frac{2u-1}{4h-1}}}\cdot \frac{1}{(n - m)^{\frac{2u-1}{4h-1}}}
\]
\[
+ \sum_{m=-K+1}^{-2K}\frac{1}{m^{1-\frac{2u}{4h-1}}} \cdot \frac{1}{(|m + K| + 1)^{\frac{2u-1}{4h-1}}}\cdot \frac{1}{(n - m)^{\frac{2u-1}{4h-1}}}
\]
\[
+ \sum_{m=-2K+1}^{n-1}\frac{1}{m^{1-\frac{2u}{4h-1}}} \cdot \frac{1}{(|m + K| + 1)^{\frac{2u-1}{4h-1}}}\cdot \frac{1}{(n - m)^{\frac{2u-1}{4h-1}}}
\]
\[
\ll \frac{1}{|K|^{\frac{2u-1}{4h-1}}} \cdot \frac{1}{n^{\frac{2u-1}{4h-1}}}\sum_{m=1}^{-\lfloor K/2 \rfloor}\frac{1}{m^{1-\frac{2u}{4h-1}}} + \frac{1}{|K|^{1-\frac{2u}{4h-1}}} \cdot \frac{1}{n^{\frac{2u-1}{4h-1}}}\sum_{m=-\lfloor K/2 \rfloor +1}^{-K}\frac{1}{(|m + K| + 1)^{\frac{2u-1}{4h-1}}}
\]
\[
+ \frac{1}{|K|^{1-\frac{2u}{4h-1}}}\cdot \frac{1}{n^{\frac{2u-1}{4h-1}}}\sum_{m=-K}^{-2K-1}\frac{1}{(|m + K| + 1)^{\frac{2u-1}{4h-1}}} + \sum_{m=-2K+1}^{n-1}\frac{1}{m^{1-\frac{2u}{4h-1}}} \cdot \frac{1}{m^{\frac{2u-1}{4h-1}}}\cdot \frac{1}{(n - m)^{\frac{2u-1}{4h-1}}}
\]
\[
\ll \frac{1}{|K|^{\frac{2u-1}{4h-1}}} \cdot \frac{1}{n^{\frac{2u-1}{4h-1}}}\int_{0}^{|K|}\frac{dx}{x^{1-\frac{2u}{4h-1}}} + \frac{1}{|K|^{1-\frac{2u}{4h-1}}} \cdot \frac{1}{n^{\frac{2u-1}{4h-1}}}\int_{0}^{|K|}\frac{dx}{x^{\frac{2u-1}{4h-1}}}
\]
\[
+ \frac{1}{|K|^{1-\frac{2u}{4h-1}}}\cdot \frac{1}{n^{\frac{2u-1}{4h-1}}}\int_{0}^{|K|}\frac{dx}{x^{\frac{2u-1}{4h-1}}} + \sum_{m=1}^{n-1}\frac{1}{m^{1-\frac{1}{4h-1}}} \cdot \frac{1}{(n - m)^{\frac{2u-1}{4h-1}}}
\]
\[
\ll \frac{|K|^{\frac{2u}{4h-1}}}{|K|^{\frac{2u-1}{4h-1}}\cdot n^{\frac{2u-1}{4h-1}}} + \frac{|K|^{1-\frac{2u-1}{4h-1}}}{|K|^{1-\frac{2u}{4h-1}}\cdot n^{\frac{2u-1}{4h-1}}} + \frac{|K|^{1-\frac{2u-1}{4h-1}}}{|K|^{1-\frac{2u}{4h-1}}\cdot n^{\frac{2u-1}{4h-1}}} + \frac{1}{n^{\frac{2u-2}{4h-1}}} = O(1).
\]

\textbf{Subase 2.3.} $K < -\frac{n}{3}$. Then we have
\[
\mathbb{E}(Y^{(\overline{\beta}, \overline{\gamma}, \overline{\delta}, \overline{\mu}, \overline{\nu})}_{u,v,h,B}(n)) \ll \sum_{m=1}^{-\lfloor K/2 \rfloor}\frac{1}{m^{1-\frac{2u}{4h-1}}} \cdot \frac{1}{(|m + K| + 1)^{\frac{2u-1}{4h-1}}}\cdot \frac{1}{(n - m)^{\frac{2u-1}{4h-1}}}
\]
\[
+ \sum_{m=-\lfloor K/2 \rfloor + 1}^{n-1}\frac{1}{m^{1-\frac{2u}{4h-1}}} \cdot \frac{1}{(|m + K| + 1)^{\frac{2u-1}{4h-1}}}\cdot \frac{1}{(n - m)^{\frac{2u-1}{4h-1}}}.
\]
If $1 \le m < -K/2$, then $|m + K| \ge |K/2| \ge \frac{n}{6}$. It follows from (i) of Lemma 4 that
\[
\sum_{m=1}^{-\lfloor K/2 \rfloor}\frac{1}{m^{1-\frac{2u}{4h-1}}} \cdot \frac{1}{(|m + K| + 1)^{\frac{2u-1}{4h-1}}}\cdot \frac{1}{(n - m)^{\frac{2u-1}{4h-1}}} \ll \sum_{m=1}^{-\lfloor K/2 \rfloor}\frac{1}{m^{1-\frac{1}{4h-1}}} \cdot \frac{1}{n^{\frac{2u-1}{4h-1}}} \cdot \frac{1}{(n - m)^{\frac{2u-1}{4h-1}}}
\]
\[
\ll \frac{1}{n^{\frac{2u-1}{4h-1}}} \cdot \sum_{m=1}^{-\lfloor K/2 \rfloor}\frac{1}{m^{1-\frac{1}{4h-1}}}\cdot \frac{1}{(n - m)^{\frac{2u-1}{4h-1}}} \ll \frac{1}{n^{\frac{2u-1}{4h-1}}\cdot n^{\frac{2u-2}{4h-1}}} = O(1).
\]

On the other hand, by $-\frac{K}{2} \ge \frac{n}{6}$,
\[
\sum_{m=-\lfloor K/2 \rfloor}^{n-1}\frac{1}{m^{1-\frac{2u}{4h-1}}} \cdot \frac{1}{(|m + K| + 1)^{\frac{2u-1}{4h-1}}}\cdot \frac{1}{(n - m)^{\frac{2u-1}{4h-1}}}
\]
\[
\ll \frac{1}{n^{1-\frac{2u}{4h-1}}} \cdot \sum_{m=-\lfloor K/2 \rfloor}^{n-1}\frac{1}{(|m + K| + 1)^{\frac{2u-1}{4h-1}}}\cdot \frac{1}{(n - m)^{\frac{2u-1}{4h-1}}}.
\]
We know that $|K|\le nh$, thus in the above sum, $1\le |m + K| + 1 \le (1 + h)n$ and $1 \le n - m\le n - 1$. Thus a positive integer among the integers $n - m$ and $|m + K| + 1$ appear, as the non-largest, at most three times. Then,
\[
\frac{1}{n^{1-\frac{2u}{4h-1}}} \cdot \sum_{m=-\lfloor K/2 \rfloor}^{n-1}\frac{1}{(|m + K| + 1)^{\frac{2u-1}{4h-1}}}\cdot \frac{1}{(n - m)^{\frac{2u-1}{4h-1}}} \le \frac{3}{n^{1-\frac{2u}{4h-1}}} \cdot \sum_{k=1}^{(h+1)n}\frac{1}{k^{\frac{2u-1}{4h-1}}}\cdot \frac{1}{k^{\frac{2u-1}{4h-1}}}
\]
\[
\ll \frac{1}{n^{1-\frac{2u}{4h-1}}} \cdot \int_{0}^{(h+1)n}\frac{dx}{x^{\frac{4u-2}{4h-1}}} \ll \frac{n^{1-\frac{4u-2}{4h-1}}}{n^{1-\frac{2u}{4h-1}}} = O(1).
\]

\textbf{Case 3.} $0 < u < v < 2h$. Then,
\[
\mathbb{E}(Y^{(\overline{\beta}, \overline{\gamma}, \overline{\delta}, \overline{\mu}, \overline{\nu})}_{u,v,h,B}(n)) \ll
\]
\begin{align*}
&\ll \sum_{\substack{(z_{1},\dots{} ,z_{2h})\in (\mathbb{Z}^{+})^{2h} \\ z_{1} + \dots{} + z_{2h} = n}}\frac{1}{(z_{1} \cdots{} z_{2h})^{\frac{4h-3}{4h-1}}}\cdot \frac{1}{(|z_{1} + \dots{} + z_{u} - (z_{u+1} + \dots{} + z_{v}) + K| + 1)^{\frac{2v-1}{4h-1}}}\\
&\ll \sum_{\substack{(m_{1},m_{2})\in (\mathbb{Z}^{+})^{2} \\ m_{1} + m_{2} < n}}\left(\sum_{\substack{(z_{1},\dots{} ,z_{u}) \in (\mathbb{Z}^{+})^{u}\\ z_{1} + \dots{} + z_{u} = m_{1}}}\left(\frac{1}{(z_{1} \cdots{} z_{u})^{\frac{4h-3}{4h-1}}}\cdot \left(\sum_{\substack{(z_{u+1},\dots{} ,z_{v})\in (\mathbb{Z}^{+})^{v-u} \\ z_{u+1} + \dots{} + z_{v} = m_{2}}}\frac{1}{(z_{u+1} \cdots{} z_{v})^{\frac{4h-3}{4h-1}}}\right)\cdot\right.\right.\\
&\left.\left.\cdot \frac{1}{(|m_{1}-m_{2} + K| + 1)^{\frac{2v-1}{4h-1}}}
\cdot\sum_{\substack{(z_{v+1},\dots{} ,z_{2h}) \in (\mathbb{Z}^{+})^{2h-v}\\ z_{v+1} + \dots{} + z_{2h} = n-m_{1}-m_{2}}}
\frac{1}{(z_{v+1} \cdots{} z_{2h})^{\frac{4h-3}{4h-1}}}\right)\right).\\
\end{align*}
It follows from (iv) of Lemma 4 that
\bigskip
\par
$\mathbb{E}(Y^{(\overline{\beta}, \overline{\gamma}, \overline{\delta}, \overline{\mu}, \overline{\nu})}_{u,v,h,B}(n)) \ll$
\begin{align*}
&\ll \sum_{\substack{(m_{1},m_{2})\in (\mathbb{Z}^{+})^{2} \\m_{1}+m_{2}<n}}\left(\sum_{\substack{(z_{1},\dots{} ,z_{u}) \in (\mathbb{Z}^{+})^{u}\\ z_{1} + \dots{} + z_{u} = m_{1}}}\left(\frac{1}{(z_{1} \cdots{} z_{u})^{\frac{4h-3}{4h-1}}}
\cdot\right.\right.\\
&\left.\left.\cdot\sum_{\substack{(z_{u+1},\dots{} ,z_{v})\in (\mathbb{Z}^{+})^{v-u} \\ z_{u+1} + \dots{} + z_{v} = m_{2}}}\frac{1}{(z_{u+1} \cdots{} z_{v})^{\frac{4h-3}{4h-1}}}\frac{1}{(|m_{1}-m_{2} + K| + 1)^{\frac{2v-1}{4h-1}}}\cdot \frac{1}{(n-m_{1}-m_{2})^{\frac{2v-1}{4h-1}}}\right)\right)\\
&\ll \sum_{\substack{(m_{1},m_{2})\in (\mathbb{Z}^{+})^{2} \\m_{1}+m_{2}<n}}\left(\sum_{\substack{(z_{1},\dots{} ,z_{u}) \in (\mathbb{Z}^{+})^{u}\\ z_{1} + \dots{} + z_{u} = m_{1}}}
\left(\frac{1}{(z_{1} \cdots{} z_{u})^{\frac{4h-3}{4h-1}}}
\cdot \frac{1}{m_{2}^{1-\frac{2(v-u)}{4h-1}}}
\cdot \frac{1}{(|m_{1}-m_{2} + K| + 1)^{\frac{2v-1}{4h-1}}}\right.\right.\\
&\left.\left.\cdot \frac{1}{(n-m_{1}-m_{2})^{\frac{2v-1}{4h-1}}}\right)\right)
\end{align*}
\[
= \sum_{\substack{(m_{1},m_{2})\in (\mathbb{Z}^{+})^{2} \\ m_{1} + m_{2} < n}}\frac{1}{m_{1}^{1-\frac{2u}{4h-1}}}\cdot \frac{1}{m_{2}^{1-\frac{2(v-u)}{4h-1}}}
\cdot \frac{1}{(|m_{1}-m_{2} + K| + 1)^{\frac{2v-1}{4h-1}}}\cdot \frac{1}{(n-m_{1}-m_{2})^{\frac{2v-1}{4h-1}}}.
\]

Now, it is enough to prove that this sum is
$O\left(\frac{n^{\frac{1}{4h-1}}}{(\log n)^{4h}}\right)$.
It follows from (i) of Lemma 4 that
\[
\sum_{\substack{(m_{1},m_{2})\in (\mathbb{Z}^{+})^{2} \\m_{1}+m_{2}<n}}\frac{1}{m_{1}^{1-\frac{2u}{4h-1}}}\cdot \frac{1}{m_{2}^{1-\frac{2(v-u)}{4h-1}}}\cdot \frac{1}{(n-m_{1}-m_{2})^{\frac{2v-1}{4h-1}}}
\]
\[
\le \sum_{m_{1}=1}^{n-1}\frac{1}{m_{1}^{1-\frac{2u}{4h-1}}}
\cdot \sum_{m_{2}=1}^{n-m_{1}-1}\frac{1}{m_{2}^{1-\frac{2(v-u)}{4h-1}}}\cdot \frac{1}{(n-m_{1}-m_{2})^{\frac{2v-1}{4h-1}}}
\]
\[
\ll \sum_{m_{1}=1}^{n-1}\frac{1}{m_{1}^{1-\frac{2u}{4h-1}}}\cdot \frac{1}{(n-m_{1})^{\frac{2u-1}{4h-1}}} \ll n^{\frac{1}{4h-1}},
\]
then
\[
\sum_{\substack{(m_{1},m_{2})\in (\mathbb{Z}^{+})^{2} \\m_{1}+m_{2}<n\\|m_{1}-m_{2} + K|>(\log n)^{16h^{2}}}}\frac{1}{m_{1}^{1-\frac{2u}{4h-1}}}\cdot \frac{1}{m_{2}^{1-\frac{2(v-u)}{4h-1}}}\cdot \frac{1}{(|m_{1}-m_{2} + K| + 1)^{\frac{2v-1}{4h-1}}}\cdot \frac{1}{(n-m_{1}-m_{2})^{\frac{2v-1}{4h-1}}}
\]
\[
\ll \frac{n^{\frac{1}{4h-1}}}{(\log n)^{4h}}.
\]
Thus it is enough to show that
\[
S = \sum_{\substack{(m_{1},m_{2})\in (\mathbb{Z}^{+})^{2} \\m_{1}+m_{2}<n\\|m_{1}-m_{2} + K|\le (\log n)^{16h^{2}}}}\frac{1}{m_{1}^{1-\frac{2u}{4h-1}}}\cdot \frac{1}{m_{2}^{1-\frac{2(v-u)}{4h-1}}}\cdot \frac{1}{(|m_{1}-m_{2} + K| + 1)^{\frac{2v-1}{4h-1}}}\cdot \frac{1}{(n-m_{1}-m_{2})^{\frac{2v-1}{4h-1}}}.
\]
\begin{equation}
\ll \frac{n^{\frac{1}{4h-1}}}{(\log n)^{4h}}
\end{equation}
We have three subcases.

\textbf{Subcase 3.1} $|K| < 10(\log n)^{64h^{3}}$. Then
\[
\sum_{\substack{(m_{1},m_{2})\in (\mathbb{Z}^{+})^{2} \\m_{1}+m_{2}<n\\|m_{1}-m_{2} + K|\le (\log n)^{16h^{2}}}}\frac{1}{m_{1}^{1-\frac{2u}{4h-1}}}\cdot \frac{1}{m_{2}^{1-\frac{2(v-u)}{4h-1}}}\cdot \frac{1}{(|m_{1}-m_{2} + K| + 1)^{\frac{2v-1}{4h-1}}}\cdot \frac{1}{(n-m_{1}-m_{2})^{\frac{2v-1}{4h-1}}}
\]
\[
= \sum_{\substack{(m_{1},m_{2})\in (\mathbb{Z}^{+})^{2} \\m_{1}+m_{2}<n\\|m_{1}-m_{2} + K|\le (\log n)^{16h^{2}}\\m_{1}<100(\log n)^{64h^{3}}}}\frac{1}{m_{1}^{1-\frac{2u}{4h-1}}}\cdot \frac{1}{m_{2}^{1-\frac{2(v-u)}{4h-1}}}\cdot \frac{1}{(|m_{1}-m_{2} + K| + 1)^{\frac{2v-1}{4h-1}}}\cdot \frac{1}{(n-m_{1}-m_{2})^{\frac{2v-1}{4h-1}}}
\]
\[
+ \sum_{\substack{(m_{1},m_{2})\in (\mathbb{Z}^{+})^{2} \\m_{1}+m_{2}<n\\|m_{1}-m_{2} + K|\le (\log n)^{16h^{2}}\\m_{1}\ge 100(\log n)^{64h^{3}}}}\frac{1}{m_{1}^{1-\frac{2u}{4h-1}}}\cdot \frac{1}{m_{2}^{1-\frac{2(v-u)}{4h-1}}}\cdot \frac{1}{(|m_{1}-m_{2} + K| + 1)^{\frac{2v-1}{4h-1}}}\cdot \frac{1}{(n-m_{1}-m_{2})^{\frac{2v-1}{4h-1}}}
\]
\[
= S_{1} + S_{2}.
\]
In $S_{1}$, there are $O((\log n)^{64h^{3}})$ possibilities for $m_{1}$, and for a fixed $m_{1}$, there are $O((\log n)^{16h^{2}})$ possibilities for $m_{2}$. Since every term in $S_{1}$ is bounded, thus
\[
S_{1} = O((\log n)^{64h^{3}+16h^{2}}) = O\left(\frac{n^{\frac{1}{4h-1}}}{(\log n)^{4h}} \right).
\]
If $m_{1}\ge 100(\log n)^{64h^{3}}$, then $|K|\le 10(\log n)^{64h^{3}}$ and $|m_{1}-m_{2} + K|\le (\log n)^{16h^{2}}$, which implies that $m_{1} \asymp m_{2}$, and so
\[
S_{2} \ll \sum_{m_{1}=1}^{n-1}\frac{1}{m_{1}^{1-\frac{2u}{4h-1}}}\cdot \frac{1}{m_{1}^{1-\frac{2(v-u)}{4h-1}}}\cdot (\log n)^{16h^{2}}
\]
because for a fixed $m_{1}$, there are $O((\log n)^{16h^{2}})$ possibilities for $m_{2}$ and
\[
\frac{1}{(|m_{1}-m_{2} + K| + 1)^{\frac{2v-1}{4h-1}}}\cdot \frac{1}{(n-m_{1}-m_{2})^{\frac{2v-1}{4h-1}}} = O(1).
\]
Then
\[
S_{2} \ll \sum_{m_{1}=1}^{\infty}\frac{1}{m_{1}^{2-\frac{2v}{4h-1}}}\cdot (\log n)^{16h^{2}} \ll (\log n)^{16h^{2}} = O\left(\frac{n^{\frac{1}{4h-1}}}{(\log n)^{4h}} \right)
\]
because $v \le 2h - 1$.

\textbf{Subase 3.2} $K \ge 10(\log n)^{64h^{3}}$. Then by $|m_{1}-m_{2} + K|\le (\log n)^{16h^{2}}$, $m_{2} \ge m_{1} + K - (\log n)^{16h^{2}} \ge m_{1}$. Thus,
\[
S_{2} \ll \sum_{m_{1}=1}^{n-1}\frac{1}{m_{1}^{1-\frac{2u}{4h-1}}}\cdot
\frac{1}{m_{1}^{1-\frac{2(v-u)}{4h-1}}}\cdot (\log n)^{16h^{2}}
= (\log n)^{16h^{2}}\cdot \sum_{m_{1}=1}^{\infty}\frac{1}{m_{1}^{2-\frac{2v}{4h-1}}}
\]
\[
= O((\log n)^{16h^{2}}) = O\left(\frac{n^{\frac{1}{4h-1}}}{(\log n)^{4h}}\right)
\]
because for a fixed $m_{1}$, there are $O((\log n)^{16h^{2}})$ possibilities for $m_{2}$ and $v \le 2h - 1$.

\textbf{Subcase 3.3} $K \le -10(\log n)^{16h^{2}}$. Then by
$|m_{1}-m_{2} + K|\le (\log n)^{16h^{2}}$,
$m_{1} \ge m_{2} - K - (\log n)^{16h^{2}} \ge m_{2}$. Then
\[
S_{2} \ll \sum_{m_{2}=1}^{n-1}\frac{1}{m_{2}^{1-\frac{2u}{4h-1}}}\cdot
\frac{1}{m_{2}^{1-\frac{2(v-u)}{4h-1}}}\cdot (\log n)^{16h^{2}} \ll (\log n)^{16h^{2}}
\cdot \sum_{m_{2}=1}^{\infty}\frac{1}{m_{2}^{2-\frac{2v}{4h-1}}}
\]
\[
= O((\log n)^{16h^{2}}) = O\left(\frac{n^{\frac{1}{4h-1}}}{(\log n)^{4h}}\right)
\]
because for a fixed $m_{2}$, there are $O((\log n)^{16h^{2}})$ possibilities for $m_{1}$ and $v \le 2h - 1$. The proof is completed.

\end{document}